%% file: 0-mesh-refinement-for-anisotropic-diffusion.tex
\documentclass[preprint,12pt]{elsarticle}

\usepackage{amssymb}
\usepackage{amsmath}
\usepackage{cleveref}
\usepackage{comment}
\usepackage[margin=1in]{geometry}
\usepackage{graphicx}
\usepackage{units}
\usepackage{xcolor}
\usepackage{marginnote}

\usepackage{lineno}

\usepackage[finalnew]{trackchanges} % use to render as resubmitted

\journal{Computers \& Mathematics with Applications}

\graphicspath{{images/}}

\newcommand{\V}[1]{\mathbf{#1}} % general vector for text
\newcommand{\Vs}[1]{\boldsymbol{#1}} % general vector for symbol
\newcommand{\Vhat}[1]{\V{\hat{#1}}} % general unit vector for text
\newcommand{\Vhats}[1]{\Vs{\hat{#1}}} % general unit vector for text

\def\density{n}				% number density
\def\Temp{T}				% temperature
\def\kpar{\kappa_{\|}} 		% parallel thermal conductivity
\def\kperp{\kappa_\perp} 		% perpendicular thermal conductivity
\def\kwedge{\kappa_\wedge} 		% perpendicular thermal conductivity
\def\kappaTensor{\Vs{\kappa}}

\def\mass{m}	%mass
\def\density{n}	%particle number density

\newcommand{\x}{\V{x}} 	% general spatial location
\def\xcoord{x}
\def\ycoord{y}
\def\zcoord{z}
\def\xhat{\hat{\V{\xcoord}}}
\def\yhat{\hat{\V{\ycoord}}}
\def\zhat{\hat{\V{\zcoord}}}

\def\Zmajor{Z}			% cylindrical coordinates
\def\Zhat{\hat{\V{\Zmajor}}}
\def\Rmajor{R}			% major radius - cylindrical coordinates
\def\aminor{a}			% largest minor radius
\def\torAngle{\zeta}		% toroidal angle
\def\pol{p}				% poloidal direction subscript
\def\tor{t}				% toroidal direction subscript

\def\Bfield{B}			% the magnetic field magnitude
\def\VBfield{\V{\Bfield}} 	% the magnetic field vector
\def\VBpol{\VBfield_\pol} 	% the poloidal magnetic field vector
\def\VBtor{\VBfield_\tor} 	% the toroidal magnetic field vector
\def\Bpol{\Bfield_\pol} 	% the poloidal component of the magnetic field
\def\Btor{\Bfield_\tor} 	% the toroidal component of the magnetic field
\def\bhat{\Vhat{b}} 		% the normalized magnetic field direction

\def\Az{A_z} 			% the z component of the vector potential

\def\hx{h_\xcoord}		% mesh spacing in x direction
\def\hy{h_\ycoord}		% mesh spacing in y direction

\def\nhat{\hat{\V{n}}} 		% unit vector normal to surface

\def\Domain{\Omega}	% Domain
\def\BCs{\mathcal{BC}}	% boundary conditions
\def\source{s}			% source
\def\dVol{d\mathcal V}	% differential volume element
\def\dArea{d\mathcal S}	% differential area element

\def\poly{p}		% polynomial order
\def\Mgrid{m}		% # grid points per dimension
\newcommand{\Hilbert}[1]{H^{#1}}  % Hilbert space of arbitrary order
\newcommand{\Sobolev}[1]{W^{#1}} % Sobolev space of arbitrary order
\def\Hone{\Hilbert{1}}	% H^1 function space
\def\Hdiv{\Hilbert{\text{div}}} % H^div finite element space
\def\Hcurl{\Hilbert{\text{curl}}} % H^curl finite element space
\def\Ltwo{L^2}
\def\fElem{E} % finite element
\def\fluxDiff{\V{e}_\sigma} % difference between numerical and analytic fluxes
\def\fluxDiffEst{\tilde{\V{e}}_\sigma} % difference between flux in Hcurl and in Hdiv

\def\uni{{\rm uniform}}
\def\isoamr{{\rm iso-amr}}
\def\anisoamr{{\rm aniso-amr}}
\def\elements{{   N}}

\def\Len{L}
\def\wrad{w}
\def\aminor{a}
\def\Rmajor{R}
\def\Rconn{R_\|}
\def\qmag{q_{mag}}

\def\dim{d}
\def\xrad{x}

\def\heatfluxVector{\mathbf q}
\def\sqrtkparkperp{\sqrt{\kpar/\kperp}}
\newcommand\norm[1]{\left| #1 \right|}

\def\constlayerwidth{\tilde{w}}

\begin{document}

\begin{frontmatter}



\title{Mesh Refinement for Anisotropic Diffusion in Magnetized Plasmas}
\author{Christopher J. Vogl\corref{cor1}\fnref{label1}}
\ead{vogl2@llnl.gov}
%\fntext[label1]{}
\cortext[cor1]{Corresponding author}
\author{Ilon Joseph\fnref{label1}}
\author{Milan Holec\fnref{label1}}

\affiliation[label1]{organization={Lawrence Livermore National Laboratory},%Department and Organization
            addressline={7000 East Ave},
            city={Livermore},
            postcode={94550},
            state={CA},
            country={U.S.}}

\begin{abstract}
\input{sections/abstract.tex}
\end{abstract}


%\item High-order finite element method discretization of anisotropic diffusion
%\item Anisotropic diffusion at realistic tokamak geometries
%\item Efficient error-based and separatrix-base adaptive mesh refinement
%\item Theory of optimal mesh refinement with respect to boundary layer


\end{frontmatter}

\nolinenumbers

\section{Introduction \label{sec:intro} }
\input{sections/introduction.tex}

\section{Steady State Plasma Transport \label{sec:plasma} }
\input{sections/plasma-transport.tex}

\section{Numerical Methods \label{sec:methods} }
\input{sections/methods.tex}

\section{Numerical Results \label{sec:results} }
\input{sections/results.tex}

\section{Discussion of Computational Efficiency \label{sec:discussion} }
\input{sections/discussion.tex}

\section{Conclusion \label{sec:conclusion} }
\input{sections/conclusion.tex}

\section*{Acknowledgements}

  The authors thank Ben Zhu, Ben Dudson, Will Pazner, and Ben Southworth for many thoughtful discussions.
  This work was performed under the auspices of the U.S. Department of Energy by Lawrence Livermore National Laboratory under Contract DE-AC52-07NA27344.
  LLNL-JRNL-840380.
  This work was supported by LLNL Laboratory Directed Research and Development project PLS-20-ERD-038.


%%       \citet{<label>} ==> Jones et al. [21]
%%       \citep{<label>} ==> [21]
%%

\bibliographystyle{elsarticle-num-names}
\bibliography{references.bib}


%\begin{thebibliography}{00}
%
%%% \bibitem[Author(year)]{label}
%%% Text of bibliographic item
%\end{thebibliography}
\end{document}

%% file: sections/abstract.tex
Highly accurate simulation of plasma transport is needed to drive the successful design and operation of magnetically confined fusion reactors.
Unfortunately, the extreme anisotropy present in magnetized plasmas results in thin boundary layers that are expensive to resolve.
This work investigates how various mesh refinement strategies might reduce that expense to allow for more efficient simulation
  \add{by comparing standard variable refinement approaches that use a field quantity to an adaptive approach that uses an error estimator}.
It is first verified that higher order discretization only realizes the proper rate of convergence once the mesh resolves the thin boundary layer, therefore motivating the focusing of refinement on the boundary layer.
  \change{Three mesh refinement strategies are investigated: one that focuses the refinement across the layer by using rectangular elements with a ratio equal to the boundary layer width, one that allows for exponential growth in mesh spacing away from the layer, and one adaptive strategy utilizing the established Zienkiewicz and Zhu error estimator.
  Across four two-dimensional test cases with high anisotropy, the adaptive mesh refinement strategy consistently achieves the same accuracy as uniform refinement using orders of magnitude less degrees of freedom.
  In the test case where the magnetic field is aligned with the mesh, the other refinement strategies also show substantial improvement in efficiency.}
  {For three two-dimensional test cases that contain characteristic features of tokamak magnetic fields, an exponential refinement strategy based on the magnetic flux function, which is the standard refinement approach in the field, is compared to an adaptive strategy utilizing the established Zienwiekicz and Zhu error estimator.
  The adaptive mesh refinement strategy consistently achieves the same accuracy using orders of magnitude less degrees of freedom than either exponential or uniform refinement.
  This result makes the adaptive refinement strategy more efficient than the exponential refinement strategy while also being more generalizable to problems with complex magnetic geometries.
  Scaling laws are derived that quantify the improvement in cost of the adaptive refinement strategy over other refinement approaches in 2D and 3D.}
  \remove{This work also includes a discussion generalizing the results to larger magnetic anisotropy ratios and to three-dimensional problems.
  It is shown that isotropic mesh refinement requires a number of degrees of freedom on the order of either the layer width (2D) or the square of the layer width (3D), whereas anisotropic refinement requires a number on the order of the log of layer width for all dimensions.
  It is also shown that the number of conjugate gradient iterations scales as a power of layer width when preconditioned with algebraic multigrid, whereas the number of iterations is independent of layer width when preconditioned with incomplete LU.
  All the results herein motivate future work in combining anisotropic adaptive mesh refinement using the Zienkiewicz and Zhu error estimator with scalable incomplete LU preconditioning strategies to reduce the total predicted computational work in problems with large anisotropy.}

%% file: sections/introduction.tex
The design of magnetically confined fusion reactors relies critically on the ability to accurately model plasma transport \cite{Wesson2011book,Boozer2005rmp}.
To this end, the magnetized confinement fusion community has developed highly sophisticated plasma transport modeling tools such as UEDGE \cite{Rognlien1999pop} and SOLPS-ITER \cite{Wiesen2015jnm}, which typically consider transport processes, including atomic and molecular physics, in a fixed background magnetic field.
Due to the interest in understanding the transport effects of perturbations to the magnetic field, general multi-fluid magnetohydrodynamics solvers, such as NIMROD \cite{Sovinec2003pop, Sovinec2004jcp}, M3D-C1 \cite{Ferraro2009jcp},
and JOREK \cite{Hoelzl2021nf} have been developed that also have the ability to model transport processes in an evolving magnetic field.
These codes are often based on finite volume and finite element approaches in order to precisely satisfy conservation laws for particles, momentum, and energy.

The physics of strongly magnetized plasmas is highly anisotropic with respect to the direction of the magnetic field.
To lowest order, a charged particle must travel in a helical orbit around a magnetic field line, which implies that plasma transport is fast along field lines, but slow across the field lines \cite{Braginskii1963rpp}.
The fastest transport process is the parallel, so-called Spitzer-H\"{a}rm  \cite{Spitzer1953pr} thermal conduction coefficient $\kpar$, due to the lightest particles, the electrons.
In contrast, the perpendicular thermal conduction coefficient $\kperp$, due to collisions is small enough that it is typically dominated by turbulent transport processes.
The ratio between the parallel and thermal conduction coefficients, $\kpar/\kperp$, can be as large as $10^9-10^{12}$ at the edge of a tokamak and gets higher still as the core temperature of 10-15 \unit{keV} is approached.
In the core, fluid transport is no longer accurate and kinetic transport processes such as neoclassical transport \cite{Hazeltine2003book} must be considered.

At the edge of a magnetized plasma, the boundary conditions conspire to generate a narrow boundary layer, traditionally called the ``scrape-off layer'' (SOL).
The width of this layer relative to the minor radius of the device scales as the square root of the inverse anisotropy ratio, $\sqrt{\kperp/\kpar}$.
Thus, even if turbulent perpendicular thermal transport reduces the anisotropy ratio to $10^6-10^8$, this still requires resolving a boundary layer width that is on the order of $10^{-4}-10^{-3}$ of the domain size.
This is very challenging for 3D simulations and can still be challenging in 2D, unless care is taken to reduce the complexity of the problem.

Within the finite element method framework, this work investigates how various approaches to refinement might capture the boundary layer dynamics while reducing computational effort away from the boundary layer (i.e., reduction in complexity of the problem).
The focus is more on mesh refinement ($h$-refinement) than on using higher-order polynomial spaces ($p$-refinement), because it is shown herein that the mesh itself must resolve the boundary layer before the expected higher-order accuracy with $p$-refinement is realized.
  \change{Both variable and adaptive mesh refinement strategies}
  {Standard variable refinement strategies using the magnetic flux function}
  \cite{Boozer2005rmp,FluxCoordinateBook}
  \add{and an adaptive strategy using an error estimation}
are leveraged with the goal of achieving the same accuracy as a uniform refinement approach but with substantially fewer degrees of freedom.
  \change{While the problems studied here use constant thermal conductivity and magnetic fields that are simple approximations of those expected in a realistic tokamak, the hope is that many of the challenges and successes found in this study will generalize to more complex nonlinear models of thermal transport and more realistic models of the magnetic field geometry.}
  {In a test case where the magnetic field is aligned with the mesh, two refinement strategies are investigated: one that focuses the refinement across the layer using rectangular elements with an aspect ratio equal to the boundary layer width and one that allows for exponential growth in mesh spacing away from the layer.
  Both strategies show substantial improvement over uniform refinement in efficiency, although the aspect-ratio strategy does not generalize to cases where the magnetic field is no longer aligned with the mesh.
  For three non-aligned, two-dimensional test cases that contain characteristic features of tokamak magnetic fields, the exponential refinement strategy based on the magnetic flux function, which is the standard refinement approach in the field, is compared to an adaptive strategy utilizing the established Zienwiekicz and Zhu error estimator.
  The adaptive mesh refinement (AMR) strategy requires orders of magnitude less degrees of freedom than either the exponential or uniform strategies.
  Thus, it is found that the AMR strategy is substantially more efficient than the standard strategies for the problems investigated.
  While those problems use constant thermal conductivity and magnetic fields that are simple approximations of those expected in a realistic tokamak, the expectation is that the ability of AMR to overcome the challenges found in this study will generalize to more complex nonlinear models of thermal transport and more realistic models of the magnetic field geometry.}

  \add{This work also includes a discussion generalizing the results to larger magnetic anisotropy ratios and to three-dimensional problems.
  It is shown that, relative to uniform refinement, isotropic AMR reduces the number of degrees of freedom (dofs) by the ratio of the layer width to domain size.
  Still, this requires a number of dofs on the order of this ratio in 2D or the square of this ratio in 3D.
  In principle, it is shown anisotropic AMR might only require a number of dofs that grows with the log of this ratio in any dimensions.
  It is also shown that the number of conjugate gradient iterations scales as a power of the ratio when preconditioned with algebraic multigrid, whereas the number of iterations is independent of the ratio when preconditioned with incomplete LU.
  All the results herein motivate the replacement of current strategies with anisotropic AMR using the Zienkiewicz and Zhu error estimator with scalable incomplete LU preconditioning strategies to reduce the total computational work in production tokamak simulation codes.}

The presentation of this work begins with background information in \Cref{sec:plasma} on steady state plasma transport in magnetically confined fusion devices, followed by the mathematic modeling of that steady state plasma transport in a variety of magnetic fields that are meant to approximate different aspects of the complex magnetic fields in a tokamak reactor, including planar and spatially varying magnetic fields.
\Cref{sec:methods} describes the finite element discretization used, including arbitrary order polynomial function spaces, as well as the mesh refinement and associate error estimation techniques used.
The resulting numerical results are presented in \Cref{sec:results}, including a quantitative measurement of computational efficiency across methods.
\Cref{sec:discussion} discusses scaling laws for degrees of freedom and condition number as means to comment on the computational cost when the anisotropy ratio is increased beyond the values in \Cref{sec:results}, as well as when three-dimensional problems are considered.
The paper concludes with a summary of the results and a brief discussion of future research directions in \Cref{sec:conclusion}.

%% file: sections/plasma-transport.tex
Plasma transport in magnetically confined fusion devices is a highly developed subject. The fluid equations in a strongly magnetized plasma are often referred to as the Braginskii equations after the influential review article \cite{Braginskii1963rpp}.
A concise introduction to these equations is given in \cite{NRLplasma}.

\subsection{Magnetic Confinement}
\begin{figure}[tbp]
    \center
    \includegraphics[width=\textwidth]{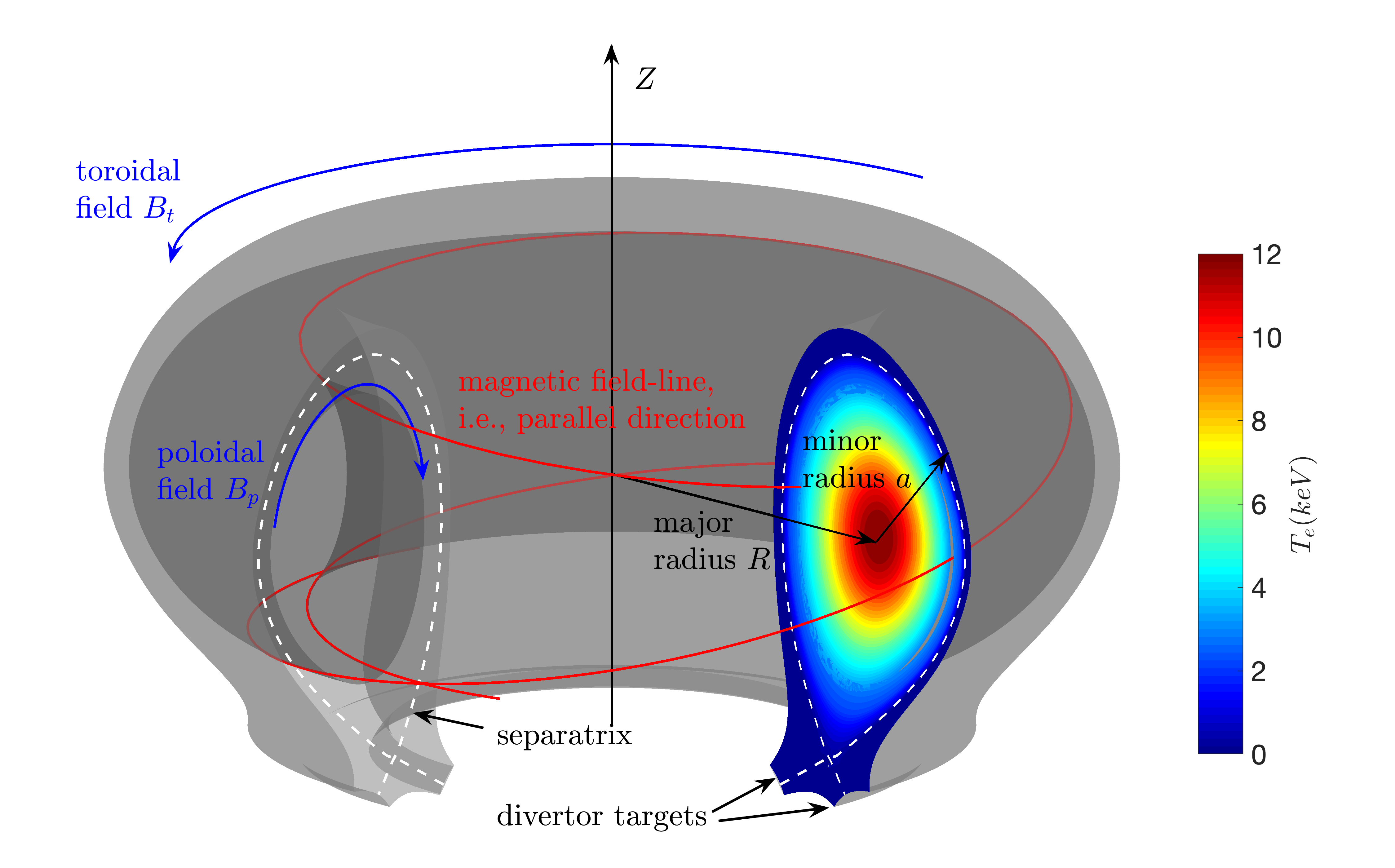}
    \caption{
      The tokamak is a toroidal magnetic confinement device with major radius $R$ and minor radius $a$. The magnetic field has a component, $B_t$, in the toroidal direction going the long way around the torus and a component, $B_p$ in the poloidal direction going the short way around the torus. A typical magnetic field line trajectory is shown in red.  The separatrix (white dashed lines) separates the region of closed field lines in the core from the region of open field lines in the scrape-off layer (SOL) that intersect the divertor target plates. The variation of the electron temperature, $T_e$, is illustrated on the cross-section of the torus to the right.
    }
    \label{fig:tokamak}
  \end{figure}

First, it is important to understand the magnetic field geometry used to confine the plasma fuel.
In a magnetic field, charged particles follow a helical gyro-motion that ensures they are tightly bound to field lines.
Bending the field lines into the shape of a torus ensures that the charged particles are trapped (to lowest order in an adiabatic expansion).
In order to control the slow drift of charged particles across field lines, as well as plasma turbulence and magnetohydrodynamic instabilities, there must be both a toroidal magnetic field, $\VBtor$, pointing the long way around the torus and a poloidal magnetic field, $\VBpol$, pointing the short way around the torus \cite{Boozer2005rmp}.
In order to have a simple yet concrete picture in mind, it is useful to consider the tokamak configuration, shown in Fig.~\ref{fig:tokamak}, which is idealized as having a perfectly axisymmetric magnetic field.
 The major radius of the torus will be denoted $\Rmajor$ and the minor radius of the torus will be denoted $\aminor$.
The rate at which field lines twist toroidally versus poloidally is called the ``safety factor'' and is approximately $\qmag\sim \Btor \aminor/\Bpol\Rmajor$.
In a typical tokamak, the ratio of $\Btor/\Bpol\sim 10$ and $\Rmajor /\aminor \sim 3$, so that the edge safety factor is of order $\qmag\sim 3-4$.

This work will focus on the important problem of predicting the temperature, $\Temp(\x)$, and heat flux, $\heatfluxVector(\x)$, within the plasma volume.
Due to the Lawson criterion, there is a range of optimal temperatures for each choice of fusion fuel \cite{Wesson2011book}.
The first generation of fusion fuels that are typically under consideration is deuterium (D) and tritium (T) precisely because they have the lowest ignition temperature, which is in the range of 10-15 \unit{keV}, as well as the highest fusion reaction rates.
Yet, the heat fluxes delivered to the surrounding walls and material surfaces must be less than 10 \unit{MW/m$^2$} in order to maintain structural integrity.
In a fusion reactor, this generally requires the regions of high plasma density near the walls to be less than 1 $\unit{eV}$ -- a factor of $10^4$ difference in temperature.
The goal of this work is to develop meshing strategies that accurately capture the variation of $T_e$ over the domain, including the core, the boundary layer outside the separatrix, and along the divertor target plates.

\subsection{Heat Flux}
For this application area, it is very important to be able to predict both the heat flux and the temperature with high accuracy.
The temperature in the interior of the plasma is important for understanding the total stored energy, but it is the heat flux near the material surfaces that is important for ensuring the reactor stays within prescribed operational limits.
The heat flux is related to the temperature gradient via the relation
\begin{gather}
	\heatfluxVector(\x) =  - \kappaTensor(\x)  \cdot \nabla \Temp(\x)
\end{gather}
where $\kappaTensor(\x)$ is the thermal conductivity tensor.
The thermal conductivity tensor has the form
\begin{gather}
	\kappaTensor(\x)=
	\kpar(\x)\bhat(\x) \otimes \bhat(\x)
	+ \kwedge (\x)\bhat(\x) \times
	+ \kperp (I - \bhat(\x) \otimes \bhat(\x))
\end{gather}
where $\bhat=\VBfield /\Bfield $ is a unit vector pointing along the magnetic field direction and $\Bfield=\norm{\VBfield}$.
According to the Braginskii equations, the components of the heat flux tensor satisfy the nonlinear scaling laws
\begin{gather}
  \begin{gathered}
	\kpar(\x)\propto \Temp^{5/2}(\x)/\mass^{1/2}\\
	\kperp (\x)\propto \density^2(\x)  \mass^{1/2}/\Bfield^2(\x) \Temp^{1/2}(\x)\\
	\kwedge  (\x)\propto \density(\x) \Temp(\x)/\Bfield(\x)
  \end{gathered}
\end{gather}
where $\mass$ is the mass of the charged particles and $\density(\x)$ is the particle number density.
The ratio of parallel to perpendicular thermal conductivity is extremely large and scales as $\kpar/\kperp \sim \Bfield^2/\mass$, which implies that it increases as the magnetic field is increased and as the mass is reduced.
Thus, this ratio is approximately 3600$\times$ smaller for electrons than for deuterium.
Because the focus of this work is on the use of advanced numerical methods for modeling plasma transport, we will eliminate the complexity associated with nonlinear variations with magnetic field, density, and temperature by treating $\kpar$ and $\kperp$ as constant in space.
More generally, the heat flux also depends on electric current, but this dependence is neglected here for simplicity.

The $\kwedge$ term in the heat flux is known as the magnetic drift term and understanding its effect is part of the subject of ``neoclassical transport'' theory \cite{HazeltineBook, HelanderBook}.
In a toroidal magnetic field, $\kwedge$ generates an effective enhancement of radial perpendicular transport that is approximately $\qmag^2 \kperp$. This represents roughly a factor of 10-16 enhancement over the $\kperp$.  For simplicity, in this work, the drift term $\kwedge$ will be neglected.

In today's fusion experiments, the actual magnitude of perpendicular thermal transport is often enhanced by orders of magnitude over the Braginskii/neoclassical levels.
This disagreement with collisional transport theory is attributed to the effect of turbulent processes which convect heat and particles from the heat source towards the edge of the plasma.
Plasma transport researchers that model the dynamic equilibrium state of a turbulent magnetized plasma experiment typically adjust the perpendicular thermal conduction coefficient $\kperp$ until experimental profiles can be matched \cite{Rognlien1999pop, Wiesen2015jnm}.
 In today's experiments, this reduces the effective anisotropy ratio to the ``measured'' range of $10^4-10^7$.

\subsection{General Test Problem}

Consider solving the linear, steady-state anisotropic diffusion problem
\begin{gather}
  \begin{gathered}
    -\nabla \cdot \big[\kappaTensor \big(\bhat(\x) \big) \cdot \nabla \Temp(\x)\big] = \source(\x), \quad \x \in \Domain, \\
    \BCs(\Temp)(\x) = 0, \quad \x \in \partial \Domain, \\
    \kappaTensor \big(\bhat(\x)\big) = \kpar \bhat(\x) \otimes \bhat(\x) + \kperp (I - \bhat(\x) \otimes \bhat(\x)),
  \end{gathered}
  \label{eq:ssp}
\end{gather}
where $\Temp$ is temperature, $\bhat$ is the magnetic field unit vector, $s$ is a source, and $\BCs$ is an operator on $\Temp$, and potentially derivatives thereof, specifying the boundary conditions.
Choices of $\bhat$, $\source$, and $\BCs$ can create boundary and internal layers of width on the order of $\sqrt{\kperp/\kpar}$, motivating an investigation into meshing strategies.

\subsection{Boundary Layer Width}

The anisotropic thermal conduction tensor in \eqref{eq:ssp} tends to generate narrow boundary layers at the edge of the plasma.
Consider, for example, the following 2D choices with $\bhat(x,y) = \xhat $ (the unit vector in the $x$ direction):
\begin{gather}
  \begin{gathered}
    -\left(\kpar \Temp_{xx}(x,y) + \kperp \Temp_{yy}(x,y)\right) = \kpar \sin(x), \quad x \in (0, \pi), \quad y \in (0,\infty),  \\
    \Temp(0,y) = \Temp(\pi,y) = \Temp(x,0) = 0, \quad \lim_{y\rightarrow \infty} \Temp_y(x,y) = 0 .
  \end{gathered}
  \label{eq:simplest}
\end{gather}
The solution to this problem is
\begin{gather}
  \Temp(x,y) = \left(1 - \exp \left[-y \sqrt{\kpar/\kperp}\right]\right) \sin(x) ,
  \label{eq:simplestSol}
\end{gather}
where the boundary layer at $y=0$ has a width proportional to $\sqrt{\kperp/\kpar}$.
Such a layer is not unique to the simplistic geometry in \eqref{eq:simplest}: similar layers are indeed expected within physical tokamaks.
%The highly anisotropic problems of interest tend to generate boundary layers that are narrow in a single radial direction, somewhat like an onion skin.
In a source-free region, the parallel and perpendicular heat fluxes must balance one another via $ \nabla \cdot  \heatfluxVector=0$.
The parallel spatial scale is set by the connection length to the wall, $\Rconn\sim \qmag \Rmajor$, while the radial boundary layer width, $\wrad$, is determined by the balance between parallel and perpendicular flows.
Assuming that the field lines are approximately straight leads to an approximately separable solution of the form $T(\x) \propto \exp{(-y/\wrad + i \ell/\Rconn)}$ where $y$ is the radial direction and $\ell$ is the length along a field line.  Thus, the balance can only hold when $\wrad/\Rconn=(\kappa_\perp/\kappa_\|)^{1/2}$.
Assume that the simulation domain length scale, $\Len$, is on the order of the parallel connection radius, $\Len \sim \Rconn $.
Then, in order to span all spatial scales from $\Len$ to $\wrad$, uniform refinement requires the number of elements
\begin{align}
  \elements_{\uni} = (\Len/\wrad)^\dim = (\kappa_\|/\kappa_\perp)^{\dim /2}.
  \label{eq:dof-scaling-uniform}
\end{align}
Thus, considering the extremely high anisotropy ratios observed for  strongly magnetized fusion plasmas, addressing the boundary layer with uniform refinement is not possible with present-day computing resources.

\subsection{Constant magnetic field problem}

  Consider $\Domain$ as the $(0,\pi) \times (0,1)$, the magnetic field $\bhat(x,y) = \Vhat{x}$, the source $\source(x,y) = \kpar \sin(x)$, and $\BCs$ such that homogeneous Dirichlet conditions are enforced in $x$ and $y$:
  \begin{gather}
    \begin{gathered}
      -\left(\kpar \Temp_{xx}(x,y) + \kperp \Temp_{yy}(x,y) \right) = \kpar \sin(x), \quad x \in (0,\pi), \quad y \in (0,1), \\
      \Temp(0,y) = \Temp(\pi,y) = \Temp(x,0) = \Temp(x,1) = 0.
    \end{gathered}
    \label{eq:toy}
  \end{gather}
  Note the exact solution to \eqref{eq:toy} is
  \begin{gather}
    \Temp(x,y) = \left(1 - \frac{\exp[-y \sqrt{\kpar/\kperp}] + \exp[-(1-y) \sqrt{\kpar/\kperp}]}{1 + \exp[-\sqrt{\kpar/\kperp}]}  \right) \sin(x).
    \label{eq:toySol}
  \end{gather}
  Defining the boundary layer width as where the corresponding exponential goes from value $1$ to $\exp[-1]$, the boundary layers at $y=0$ and $y=1$ have width $\sqrt{\kperp/\kpar}$.
  For scenarios where it may be more natural to think of forcing on $\partial \Domain$ rather than in $\Domain$, consider that $\tilde{\Temp}(x,y) = \Temp(x,y) - \sin(x)$ satisfies
  \begin{gather*}
    -\left(\kpar \tilde{\Temp}_{xx}(x,y) + \kperp \tilde{\Temp}_{yy}(x,y) \right) = 0, \quad x \in (0,\pi), \quad y \in (0,1), \\
    \tilde{\Temp}(0,y) = \tilde{\Temp}(\pi,y) = 0, \quad \tilde{\Temp}(x,0) = \tilde{\Temp}(x,1) = -\sin(x).
  \end{gather*}

\subsection{Poloidally varying magnetic field problems}

  While the constant magnetic field in \eqref{eq:toy} leads to an exact solution, it is indeed very simplistic when compared to fields expected in a magnetically confined fusion reactor.
  Recall that the magnetic field $\VBfield$ in such a device must be topologically toroidal and that the total magnetic field is the combination of the poloidal and toroidal fields: $\VBfield(\x) = \VBpol(\x)  + \VBtor(\x)$.
  Let us introduce the cylindrical coordinates $\x = \Rmajor \cos(\torAngle) \xhat+ \Rmajor \sin(\torAngle) \yhat  +\Zmajor \Zhat$, where $\torAngle$ is the toroidal angle.
  The toroidal field induced by external coils is much larger than that generated by plasma currents, so the toroidal field is approximately $\VBtor \simeq B_{t0}r_0 \Vhats{\torAngle} / \Rmajor$ where $ B_{t0}r_0$ is constant.
  The poloidal field is then defined via $\VBpol= \nabla \times \psi(\x)  \nabla \torAngle$ where $\psi(\x)$ is the poloidal magnetic function; i.e. it measures the magnetic flux through a poloidal section of the torus.
  The ideal MHD equilibrium conditions \cite{Boozer2005rmp}, which require nested toroidal flux surfaces, imply that the field lines must lie in surfaces of constant $\psi(\x)$.

\begin{figure}[tbp]
  \centering
    \includegraphics[width=2in]{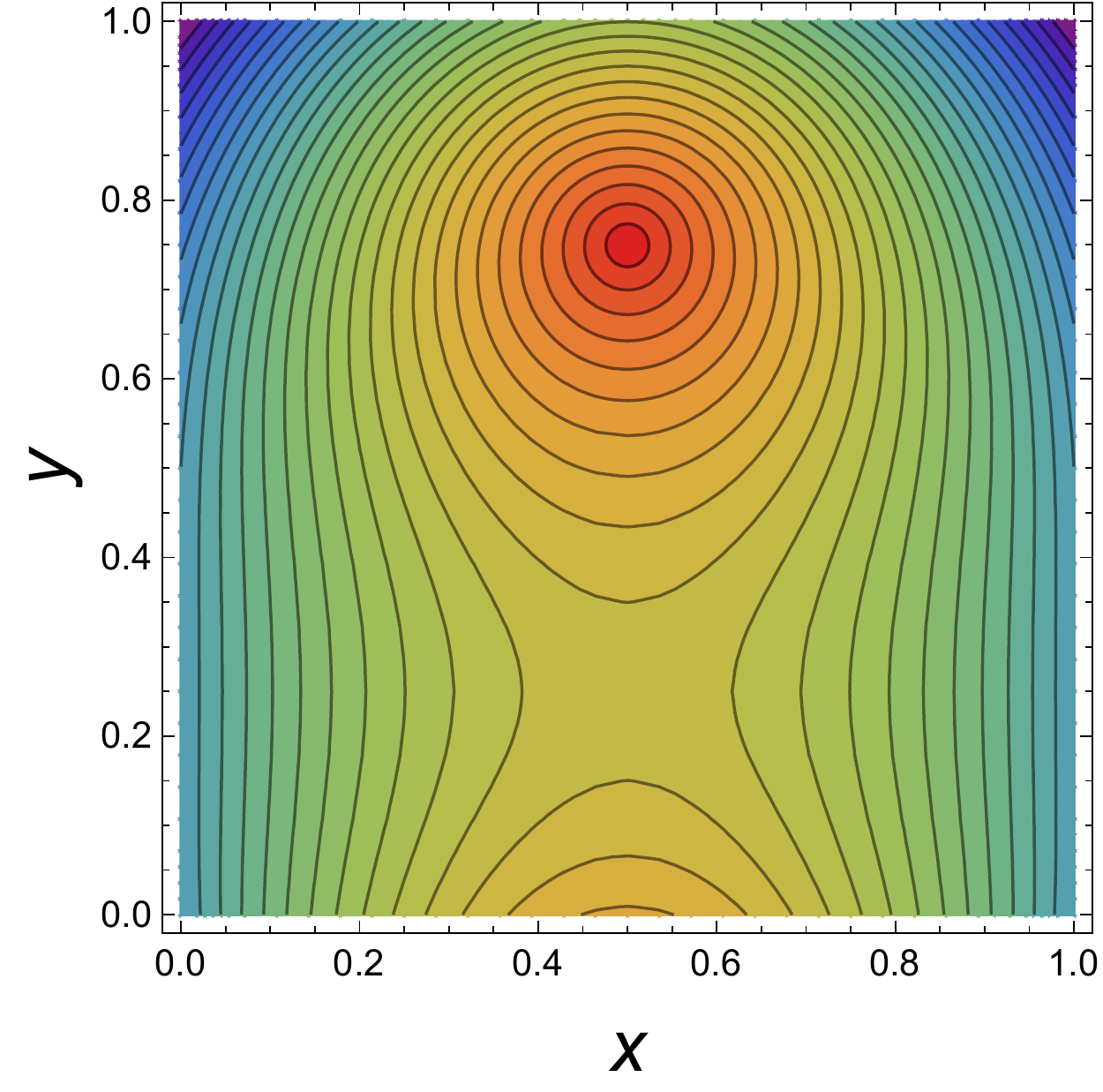}
        \includegraphics[width=2in]{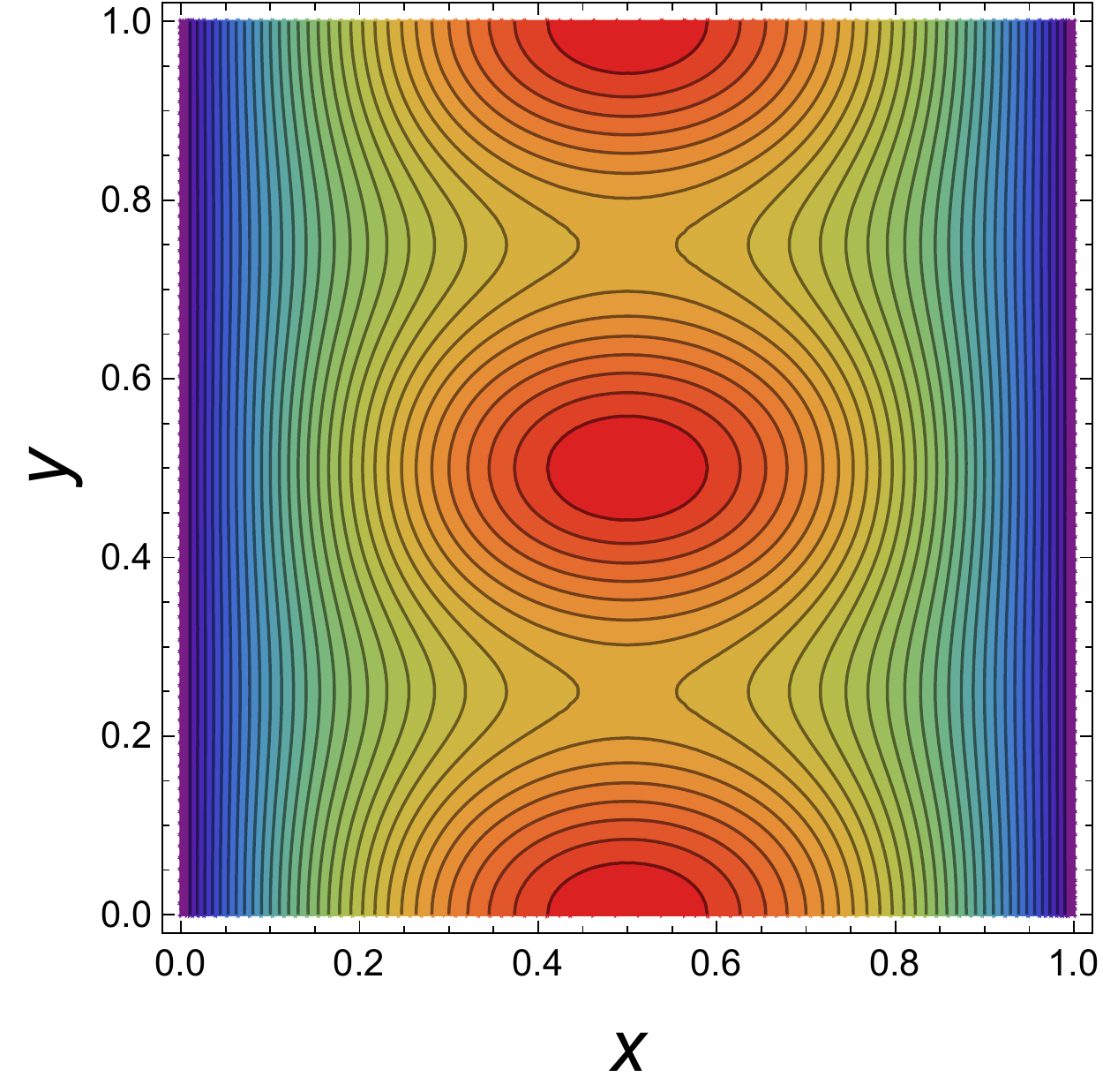}
        \includegraphics[width=2in]{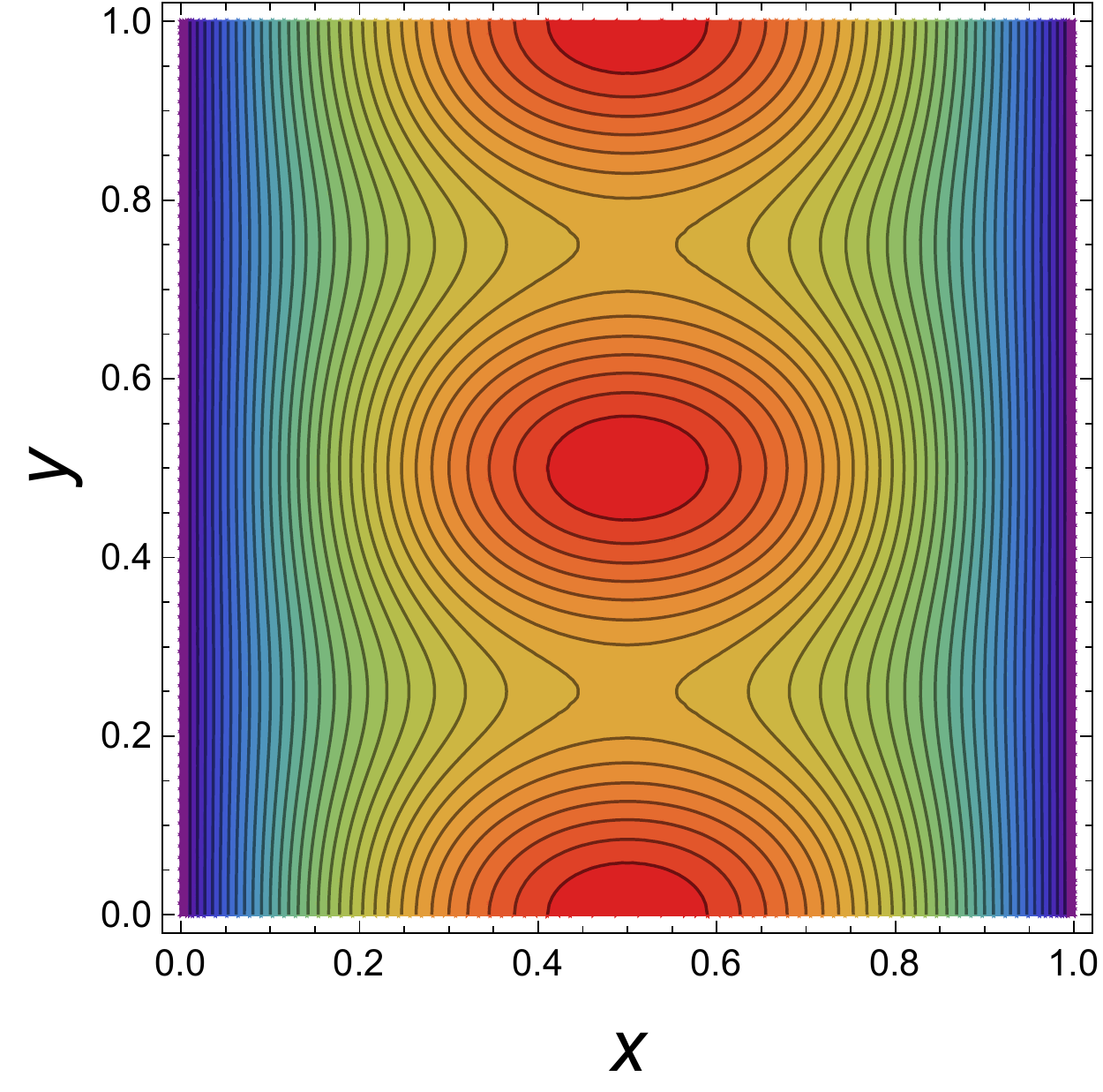}\\
     \hspace{0.25in}   (a) \hspace{1.75in} (b) \hspace{1.75in} (c)
    \caption{
Contour plots showing magnetic flux surfaces for poloidally varying problems, specifically (a) $-\exp(A_z)$ for single-null diverted tokamak, (b) $-A_z$ for double-null diverted tokamak, and (c) $-A_z$ for magnetic island.
    }
    \label{fig:magnetic flux surfaces}
  \end{figure}

  In order to further simplify the problem to focus on the key issues, we will assume a large aspect ratio toroidal geometry, where the minor radius is much smaller than the major radius $\Rmajor$; i.e. so that the major radius can be assumed to be approximately constant.
  In this case, one can assume that the toroidal field is constant $\VBtor = \zhat$ and that the poloidal field is given by $\VBpol = \nabla \times [\Az(x,y) \Vhat{z}]$, where $\psi\simeq \Az\Rmajor$.
  Furthermore, we will assume that variations along the toroidal direction now given by $\zhat$ vanish, so that
   $\bhat(x,y)\cdot \nabla = \VBpol(x,y)/|\VBfield(x,y)|\cdot \nabla $.
   Three such magnetic fields, illustrated in Fig.~\ref{fig:magnetic flux surfaces}, are now defined along with the associated test problem to be investigated.

  \subsubsection{Single-null magnetic field problem}
    The first poloidally varying magnetic field we shall study is the ``two-wire model'' of  a single null diverted tokamak, with the magnetic flux function
    \begin{gather}
      \Az(\x) = \log(|\x-\x_1||\x-\x_2|).
      \label{eq:two-wire-field}
    \end{gather}
    Note the point $\x_s = \tfrac{1}{2}(\x_1 + \x_2)$ represents the null point where the poloidal magnetic field vanishes.
    The level curve through $\x_s$, called the separatrix, partitions $\Domain$ into level curves that form closed field lines in the core that do not intersect the boundary and open field lines at the edge that do reach the boundary.
    Generally, the level curves consist of points $\x$ where the product of distances from $\x$ to $\x_1$ and to $\x_2$ remains constant.
    Considering $\Domain$ to be the unit square, the specific problem to be investigated is
    \begin{gather}
      \begin{gathered}
        -\nabla \cdot \big[ \kappaTensor\big(\bhat(x,y)\big) \cdot \nabla \Temp(x,y)\big] = e^{-\tfrac{1}{2}\big(\left(\frac{x-1/2}{1/8}\right)^2 + \left( \frac{y-1/2}{1/8}\right)^2\big)}, \quad x \in (0,1), \quad y \in (0,1), \\
        \Temp(0,y) = \Temp(1,y) = \Temp(x,0) = \Temp(x,1) = 0,
      \end{gathered}
      \label{eq:two-wire}
    \end{gather}
    Placing a heat source,  $\source(x,y)$,  inside the separatrix, traps the heat in the middle of the domain, because it can only slowly diffuse outward at the much smaller rate set by $\kappa_\perp$.
    Note that the magnetic field geometry drives an internal temperature layer of width on the order of $\sqrt{\kperp/\kpar}$ near the separatrix.

  \subsubsection{Double null magnetic field problem}
    The second poloidally varying magnetic field we shall study models the field in a double null diverted tokamak, with the magnetic flux function
    \begin{gather}
      \Az(x,y) = \tfrac{1}{2} (x - x_0)^2 + \tfrac{1}{2} \left(\tfrac{1}{4} \sin\big[2\pi (y-y_0)\big]\right)^2.
      \label{eq:dbl-null-field}
    \end{gather}
    Note the points $\x_s^\pm = (x_0,y_0 \pm \tfrac{1}{4})$ act as the two poloidal field nulls.  The  separatrix curve that passes through those points partitions $\Domain$ in a similar fashion to the single null in \eqref{eq:two-wire-field}.
    Again considering $\Domain$ to be the unit square, the specific problem to be investigated is
    \begin{gather}
      \begin{gathered}
        -\nabla \cdot \big[ \kappaTensor\big(\bhat(x,y)\big) \cdot \nabla \Temp(x,y)\big] = e^{-\tfrac{1}{2}\big(\left(\frac{x-1/2}{1/8}\right)^2 + \left(\frac{y-1/2}{1/8}\right)^2\big)}, \quad x \in (0,1), \quad y \in (0,1), \\
        \Temp(0,y) = \Temp(1,y) = \Temp(x,0) = \Temp(x,1) = 0,\\
      \end{gathered}
      \label{eq:dbl-null}
    \end{gather}
    and $\Az(x,y)$ defined by \eqref{eq:dbl-null-field} with $(x_0,y_0) = (\tfrac{1}{2},\tfrac{1}{2})$.
    As with \eqref{eq:two-wire}, placing a source inside the separatrix focuses the heat generated by $\source(x,y)$ in the middle of the domain.
    The difference between this problem \eqref{eq:dbl-null} and \eqref{eq:two-wire} is that there are two nulls near which the poloidal magnetic field vanishes. Again, internal boundary layers of width on the order of $\sqrt{\kperp/\kpar}$ are generated near the separatrix.

  \subsubsection{Magnetic island problem}
    The third poloidally varying magnetic field we will study models the way that a perturbation to the magnetic field can cause a topological change in the flux surfaces that produces a ``magnetic island.''
    The magnetic flux function to model such a phenomenon is a modified version of \eqref{eq:dbl-null-field}, where the perturbation vanishes at the radial boundaries, set by $x=0$ and $x=L$:
    \begin{gather}
      \Az(x,y) = \tfrac{1}{2} (x - x_0)^2 + \tfrac{1}{2} \left ( (1 - \tfrac{x-x_0}{L/2})(1 + \tfrac{x-x_0}{L/2}) \tfrac{1}{4} \sin \big[ 2\pi (y-y_0)\big] \right)^2.
      \label{eq:mag-island-field}
    \end{gather}
    The magnetic island width is defined by the separatrix which passes through the points $\x_s^\pm = (x_0,y_0 \pm \tfrac{1}{4})$, formed on the background magnetic flux function $\Az(x,y) = \tfrac{1}{2} (x -x_0)^2$.
    Again considering $\Domain$ to be the unit square, the specific problem to be investigated is
    \begin{gather}
      \begin{gathered}
        -\nabla \cdot \big[ \kappaTensor\big(\bhat(x,y)\big) \cdot \nabla \Temp(x,y)\big] = 0, \quad x \in (0,1), \quad y \in (0,1), \\
        \Temp(0,y) = 1, \quad \Temp(1,y) = 0, \\
        \Temp(x,0) = \Temp(x,1), \quad
        \Temp_y(x,0) = \Temp_y(x,1), \\
      \end{gathered}
      \label{eq:mag-island}
    \end{gather}
    and $\Az(x,y)$ defined by \eqref{eq:mag-island-field} with $(x_0,y_0) = (\tfrac{1}{2},\tfrac{1}{2})$ and $L=1$.
    For this problem, there is an effective heat source at the $x=0$ boundary, which maintains the temperature differential, and heat must pass through the magnetic island to reach the effective sink at $x=1$.

%% file: sections/methods.tex
The general steady state heat problem \eqref{eq:ssp} is discretized in space using the finite element method.
The corresponding weak/variational form of \eqref{eq:ssp} is that
\begin{gather}
  \begin{gathered}
    \int_\Domain \left[\kappaTensor(\bhat) \cdot \nabla \Temp \right] \cdot \nabla \psi \,\dVol - \int_{\partial \Domain} \psi \left[\kappaTensor(\bhat) \cdot \nabla \Temp \right] \cdot \nhat \,\dArea = \int_\Domain \source \psi \,\dVol
    \quad \forall \psi(\x) \in \Hone(\Domain), \\
    \Temp(\x) \in \mathcal{V}(\Domain) = \{f(\x) \in \Hone(\Domain) : \BCs(f)(\x) \text{ for } \x \in \partial \Domain  = 0\}
  \end{gathered}
  \label{eq:ssp-weak}
\end{gather}
where $\Hone(\Domain)$ is a Hilbert space.
Denote $\Domain_h$ as a collection of elements (mesh) that discretizes $\Domain$.
The space of polynomials of order $\poly$ that are piecewise defined on $\Domain_h$ is denoted $\Hone_\poly(\Domain_h)$.
Denote $\mathcal{V}_h^\poly(\Domain_h) = \{f(\x) \in \Hone_\poly(\Domain_h) : \BCs(f)(\x) = 0$ for $\x \in \partial \Domain \}$ and note that $\mathcal{V}_h^\poly(\Domain_h) \subset \mathcal{V}(\Domain)$.
The exact solution to \eqref{eq:ssp-weak} is thus approximated by $\Temp(\x) \in \mathcal{V}_h^p(\Domain)$: $\Temp(\x) = \sum_{i=1}^n \Temp_i \psi_i(\x)$ where $\{\psi_1(\x),\ldots,\psi_n(\x)\}$ is a set of basis functions for $\mathcal{V}_h^\poly(\Domain_h)$.
Because $\mathcal{V}_h^\poly(\Domain_h)$ is also a subspace of $\Hone(\Domain)$, the governing equations, $i = 1, \ldots, n$, for the approximate solution $\Temp(\x)$ are chosen to be
\begin{gather}
  \sum_{j=1}^n \left[\int_\Domain \nabla \psi_i \cdot \left[\kappaTensor(\bhat) \cdot \nabla \psi_j\right] \,\dVol - \int_{\partial \Domain} \psi_i \left[\kappaTensor(\bhat) \cdot \nabla \psi_j\right] \cdot \nhat \,\dArea \right] \Temp_j= \int_{\Domain} \source \psi_i \,\dVol .
  \label{eq:ssp-fem}
\end{gather}
The finite element problem \eqref{eq:ssp-fem} is assembled and solved using the MFEM software library \cite{mfem}.
MFEM provides a variety of arbitrary-order function spaces for unstructured meshes, including $\Hone_\poly(\Domain_h)$. The conjugate gradient method implemented in MFEM is used to solve the system, preconditioned with BoomerAMG using the provided interface to hypre \cite{falgout_hypre_2002}.

\subsection{Error Estimation}

  The finite element method provides a variety of advantages that include error estimation techniques.
  One such estimation technique for the steady state problem \eqref{eq:ssp} that utilizes various function spaces commonly used by the finite element method is from the work of Zienkiewicz and Zhu \cite{zienkiewicz_superconvergent_1992a, zienkiewicz_superconvergent_1992b, zienkiewicz_simple_1987}, commonly referred to as ZZ estimation.
  The key observation is that the flux $-\kappaTensor \nabla \Temp$ in \eqref{eq:ssp} should be continuous across element edges/faces to make the weak divergence of the flux on an element $\fElem$, evaluated as $\int_{\partial \fElem} \kappaTensor \nabla \Temp \cdot \nhat \,\dArea$, well-defined.
  In other words, the flux $-\kappaTensor \nabla \Temp$ should be in $\Hdiv(\Domain_h)$.
  The temperature field $\Temp$ is typically chosen to be in $\Hone_\poly$, which results in $\nabla \Temp \in \Hcurl_\poly(\Domain_h)$.

  Zienkiewicz and Zhu surmised that the difference $\fluxDiffEst$ between $\kappaTensor \nabla \Temp$ and its projection to $\Hdiv(\Domain_h)$ will provide an estimate for the difference $\fluxDiff$ between the numerical and analytic fluxes.
  Specifically for $\Temp \in \Hone_\poly$, the value of $\fluxDiffEst$ is defined as $\kappaTensor \nabla \Temp - \sum_{i=1}^{\tilde{n}} \sigma_i \boldsymbol{\psi}_i$, where $\{\boldsymbol{\psi}_1, \ldots, \boldsymbol{\psi}_{\tilde{n}}\}$ is a set of basis functions for $\Hdiv_{\poly-1}(\Domain_h)$ and $\{\sigma_1,\ldots,\sigma_{\tilde{n}}\}$ is determined by the projection of $\kappaTensor \nabla \Temp$ on $\Hdiv_{\poly-1}(\Domain_h)$.
  The MFEM implementation of the ZZ error estimate uses an $L_2$ projection defined by $\{\sigma_1,\ldots,\sigma_{\tilde{n}}\}$ satisfying $\int_\Omega \kappaTensor \nabla \Temp \cdot \boldsymbol{\psi}_i \, \dVol = \sum_{j=1}^{\tilde{n}} \sigma_j \boldsymbol{\psi}_j \cdot \boldsymbol{\psi}_i$ for $i=1,\ldots,\tilde{n}$.
  The error estimate is then computed per element $\fElem$ as
  \begin{gather}
    | \fluxDiffEst |_\fElem = \int_\fElem \| \fluxDiffEst \|_2 \, \dVol
    \label{eq:flux-error}
  \end{gather}
  It is worth noting if $\BCs$ enforces Dirichlet or Neumann conditions, then the energy norm of the error in $\Temp$ is equivalent to $-\int_\Domain \fluxDiff \cdot \kappaTensor^{-1} \fluxDiff \, \dVol$, and therefore can be approximated by $\fluxDiffEst$.

\subsection{Mesh Refinement}
\label{sec:methods-amr}

  MFEM supports non-conforming refinement of unstructured meshes given a refinement strategy that identifies elements for refinement or for de-refinement/coarsening.
  The refinement strategy can depend on the current solution (adaptive refinement) or it can depend on some \textit{a priori} known variables (variable refinement).
  This work utilizes one adaptive strategy and one variable strategy for identifying elements for refinement.
  The adaptive refinement strategy marks elements where the current ZZ error estimate for flux error \eqref{eq:flux-error} is greater than $\tfrac{3}{4} \tfrac{1}{N}\sum_{i=1}^N |\fluxDiffEst|_{\fElem_i}$, $\{\fElem_1,\ldots,\fElem_N\}$ being the collection of elements that define $\Domain_h$.
  As such, this strategy will be referred to as ZZ refinement.

  The variable strategy marks elements that either (i) contain the separatrix curve or (ii) have not reached a target size relative to that of the elements containing the separatrix.
  For (ii), the ratio of the target size and the size of elements containing the separatrix is given by an exponential function that grows proportional to some measure of distance of the element from the separatrix.
  As such, this strategy will be referred to as exponential refinement.
  The exponential growth is motivated by the exact solution \eqref{eq:simplestSol} and will be derived from approximation theory in Section~\ref{sec:results}.
  Note that the variable refinement strategy does require the location of the separatrix to be known \textit{a priori}, making it the less general of the two strategies.

%% file: sections/results.tex
The finite-element approximation \eqref{eq:ssp-fem} to the steady-state anisotropic diffusion problem \eqref{eq:ssp} is now solved with various magnetic fields to investigate the ability of mesh refinement to improve the computational efficiency for large $\kpar/\kperp$ values.
First, the implementation of the discretization \eqref{eq:ssp-fem} is verified on the  constant magnetic field problem \eqref{eq:toy} with convergence testing against the analytic solution \eqref{eq:toySol}.
The efficiency of using a mesh whose elements have an aspect ratio equal to the magnetic anisotropy and of using exponential refinement is investigated for the field-aligned mesh in \eqref{eq:toy}.
Finally, for the non-field-aligned single null magnetic field \eqref{eq:two-wire}, double null magnetic field \eqref{eq:dbl-null}, and magnetic island \eqref{eq:mag-island} problems, the efficiency of both ZZ and exponential refinement is compared to that of uniform refinement.

\subsection{Verification}
\label{sec:verification}

\begin{figure}[tbp]
    \center
    \includegraphics[width=\textwidth]{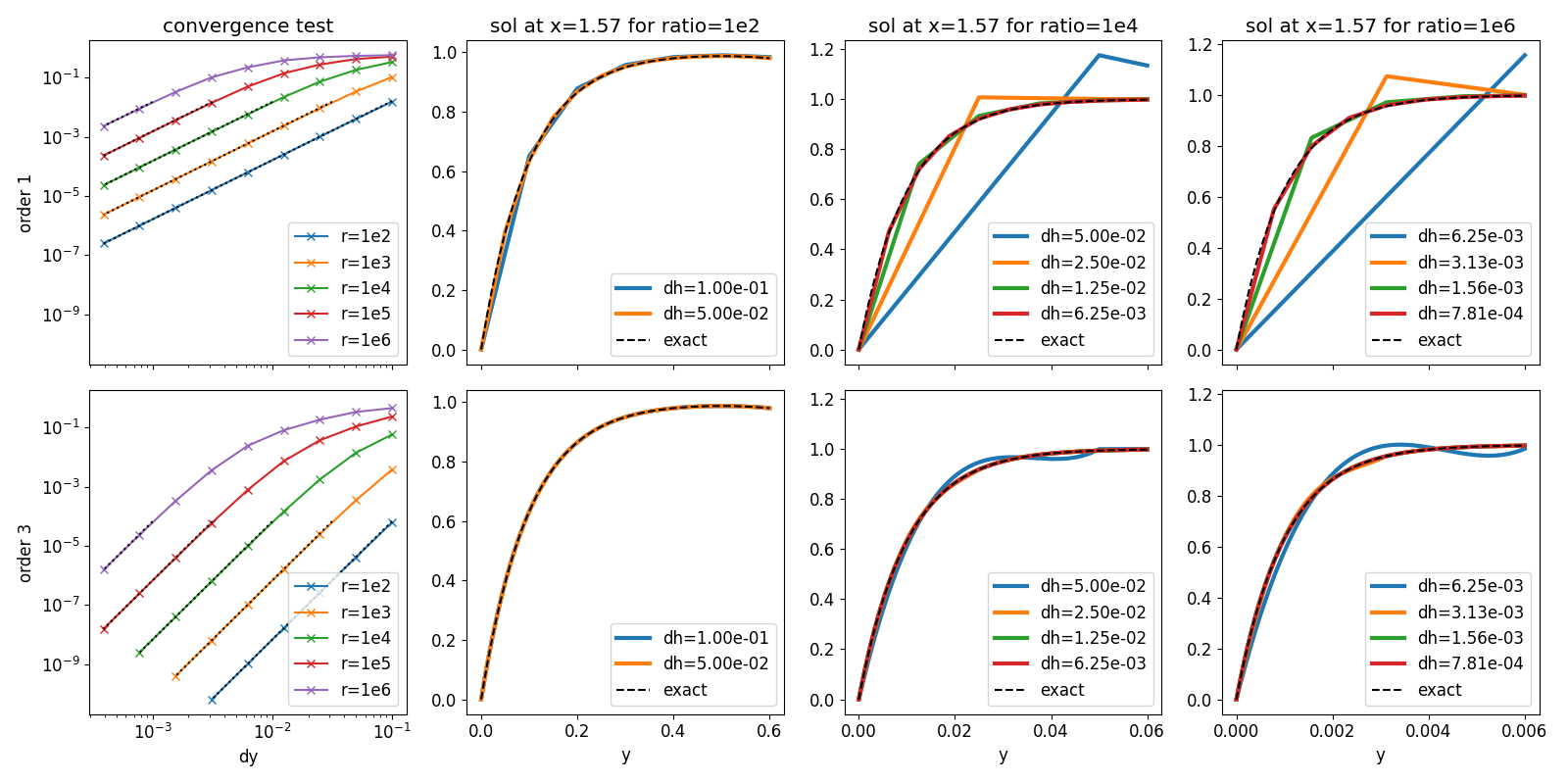}
    \caption{
      Convergence test and corresponding solution profiles for constant magnetic field problem \eqref{eq:toy} with various $\kpar/\kperp$ ratios (denoted in figure as $r$).
      Convergence test shows that expected convergence rates are not attained for either piecewise linear ($\poly=1$) and cubic ($\poly=3$) function spaces until mesh elements are smaller than the boundary layer width.
      Solution profiles show over/under-shooting behavior in both function spaces until the mesh elements are smaller than the boundary layer width, with the over/under-shooting in the $p=3$ solution being more subdued.}
    \label{fig:convergence-test}
  \end{figure}

  The implementation of the discretization \eqref{eq:ssp-fem} is first verified by a convergence test on the constant magnetic field problem \eqref{eq:toy}.
  The mesh $\Domain_h$ partitions $\Domain$ into $\Mgrid^2$ uniform, rectangular elements (i.e. that have edges of length $\hx=\pi/\Mgrid$ and $\hy=1/\Mgrid$).
  The resulting linear systems are solved with conjugate gradient preconditioned by hypre's BoomerAMG across $6$ MPI ranks.
  The convergence test uses an error norm that focuses on the discrete $L^2$ error in and around the boundary layer at $y=0$:
  \begin{gather*}
    \| \Temp - \Temp^*\| = \sqrt{ \frac{1}{15000} \sum_{i=0}^{100} \sum_{j=0}^{150} \big( \Temp(x_i, y_j) - \Temp^*(x_i,y_j)\big)^2}, \quad x_i = \frac{i}{100} \pi,\quad y_j = \frac{j}{50} \sqrt{\frac{\kperp}{\kpar}},
  \end{gather*}
  where $\Temp^*(x,y)$ is the exact solution to \eqref{eq:ssp-fem}.
  The test is conducted with piecewise linear ($\poly=1$) and piecewise cubic ($\poly=3$) function spaces for $\Mgrid=10,20,\ldots,2560$ and $\kpar/\kperp=10^2,10^3,\ldots,10^6$.
  Figure~\ref{fig:convergence-test} verifies that both the linear and cubic approximations of \eqref{eq:toySol} reach the theoretically predicted convergence order once $\hx=\hy$ is within the asymptotic regime.
  Note that the asymptotic regimes in Figure~\ref{fig:convergence-test} are consistent with requiring the mesh resolve the boundary layers, that $\hy \leq \sqrt{\kperp/\kpar}$, regardless of whether $\poly=1$ or $\poly=3$.
  The solution profiles in Figure~\ref{fig:convergence-test} also indicate that the mesh itself must resolve the boundary layer to avoid overshoot behavior in the numerical solution.

  Approximation theory can explain the over/under-shooting behavior and why one might expect the requirement that the mesh must resolve the boundary layer before the appropriate convergence rate is obtained.
  Consider the result from Ciarlet and Raviart \cite{ciarlet_general_1972} that provides the bound
  \begin{gather}
    \begin{gathered}
      \| u - \Pi_\poly u \|_{\Sobolev{k,l}(E)} \leq C h^{\poly+1-k} \left[ \int_E \sum_{a+b=p+1} \left| \frac{\partial^{\poly+1} u}{\partial x^a \partial y^b} \right|^l \dVol \right]^{1/l}, \\
      \text{where }
      \| v \|_{\Sobolev{k,l}(E)} := \left[ \int_E \sum_{c=0}^k \sum_{a+b=c} \left| \frac{\partial^c v}{\partial x^a \partial y^b} \right|^l \dVol \right]^{1/l} \quad \forall v \in \Sobolev{k,l}(E).
    \end{gathered}
    \label{eq:interpolation-error}
  \end{gather}
  for any element $E$ in $\Domain_h$, where $u$ is a function in Sobolev space $\Sobolev{\poly+1,l}(\Domain)$, $\Pi_\poly u$ is the piecewise defined polynomial of order $\poly$ that interpolates $u$.
  Because the exact solutions \eqref{eq:simplestSol} and \eqref{eq:toySol} are in $\Hilbert{\poly+1}(\Domain)$, the bound applies with $l=2$, $u=T$, and $\Pi_\poly T \in \mathcal{V}_h^\poly(\Domain_h)$.
  Because the solution \eqref{eq:simplestSol} lends itself to more compact analysis with results that can generalize to \eqref{eq:toySol} that has similar boundary layer structure, consider $\kpar/\kperp \gg 1$ for $T$ from \eqref{eq:simplestSol}:
  \begin{gather*}
    \left[ \int_E \sum_{a+b=p+1} \left | \frac{\partial^{\poly+1} \Temp}{\partial x^a \partial y^b} \right|^2 \dVol \right]^{1/2}
      \sim \left( \sqrt{\frac{\kpar}{\kperp}} \right)^{p+1} \left[ \int_E \left |\exp[-y \sqrt{\kpar/\kperp}]\sin(x) \right |^2 \right]^{1/2} \\
    \Rightarrow h^{p+1} \left[ \int_E \sum_{a+b=p+1} \left | \frac{\partial^{\poly+1} \Temp}{\partial x^a \partial y^b} \right|^2 \dVol \right]^{1/2}
      \sim \left( h \sqrt{\frac{\kpar}{\kperp}} \right)^{p+1} \exp[-\tilde{y}\sqrt{\kpar/\kperp}] \left[ \int_E \left |\sin(x) \right |^2 \dVol \right]^{1/2},
  \end{gather*}
  for some $\tilde{y}$ in element $E$ per the integral Mean Value Theorem.
  Thus, for the discrete $\Ltwo$ error (i.e., the $\Sobolev{0,2}(\Omega_h)$ error) to asymptotically decrease at the appropriate rate, one must have $h \sqrtkparkperp < 1$.

\subsection{Mesh refinement for constant magnetic field problem}

  As seen above, the requirement of the mesh resolving the boundary layer to (i) avoid overshoot behavior in the solution and (ii) capitalize on convergence rates of higher-order polynomial spaces becomes very expensive for uniform meshes.
  Consider that piecewise polynomial function spaces of order $\poly$ require $\poly+1$ degrees of freedom per element in each direction.
  Continuous function spaces, such as $\Hone_{\poly}$, share degrees of freedom across elements so that $(\poly/\hx+1)(\poly/\hy+1)$ degrees of freedom are required for $\Omega_h$ containing $(1/\hx)(1/\hy)$ elements.
  Thus, resolving the boundary layer of width $\sqrt{\kperp/\kpar}$, denoted here as $\constlayerwidth$, requires on the order of $\kpar/\kperp$ degrees of freedom.
  For field-aligned problems, one might consider using meshes that are more refined across the boundary layer than along the layer.
  Such meshes are studied for the constant magnetic field problem \eqref{eq:toy} using $\kpar/\kperp = 10^4$ and a piecewise cubic function space.
  To start, instead of the rectangles in the uniform mesh, rectangles that have an aspect ratio closer to the anisotropy ratio are used.
  Thus, $\Domain_h$ partitions $\Domain$ into $m_x \times m_y$ rectangles with $m_x$ rectangles along the $x$ direction and $m_y$ rectangles along the $y$ direction with $m_y/m_x = 100/3$ (i.e., $h_x/h_y = 100 \pi/3 \approx 100$).
  A second approach builds on the first by allowing $h_y$ to vary with $y$.
  The form of $h_y(y)$ is motivated by an observation using the approximation theory bound \eqref{eq:interpolation-error} for $\tilde{T}(y) = T(x^*,y)$ with $x^* \in (0,\pi)$ and $T(x,y)$ from \eqref{eq:simplestSol}.
  Consider the one-dimensional element $E = (y^E,y^E+h_y^E)$ where
  \begin{gather*}
    \begin{aligned}
      \left[ \int_{y^E}^{y^E+h_y^E} \left |\frac{d^{p+1}\tilde{T}}{dy^{p+1}}\right |^2 dy \right]^{1/2}
        &= \left(\frac{1}{\constlayerwidth}\right)^{p+1} \left[ \int_{y^E}^{y^E+h_y^E} \left |\exp[-y / \constlayerwidth] \sin(x^*) \right |^2 dy \right]^{1/2} \\
        &\leq \left(\frac{1}{\constlayerwidth}\right)^{p+1} \left[ \int_{y^E}^{y^E+h_y^E} \left |\exp [-y / \constlayerwidth] \right |^2 dy \right]^{1/2} \\
        &= \left(\frac{1}{\constlayerwidth}\right)^{p+1} \exp[-y^E/\constlayerwidth] \left[ \frac{1 - \exp[-2h_y^E/\constlayerwidth]}{2 h_y^E / \constlayerwidth} \int_{y^E}^{y^E+h_y^E} dy \right]^{1/2}
    \end{aligned} \\
    \Rightarrow
      \|\tilde{T} - \Pi_p \tilde{T} \|_{\Ltwo(E)} \leq C \left(\frac{h_y^E}{\constlayerwidth} \right)^{p+1} \exp[-y^E / \constlayerwidth] \left[ \frac{1 - \exp[-2h_y^E / \constlayerwidth]}{2 h_y^E / \constlayerwidth} \int_{y^E}^{y^E+h_y^E} dy \right]^{1/2} \\
    \Rightarrow
      \|\tilde{T} - \Pi_p \tilde{T} \|_{\Ltwo(\Domain_h)}^2 \leq
      \sum_{E \in \Domain_h} \left[ C \left(\frac{h_y^E}{\constlayerwidth}\right)^{p+1} \exp[-y^E / \constlayerwidth] \right]^2 \frac{1 - \exp[-2h_y^E / \constlayerwidth]}{2 h_y^E / \constlayerwidth} \int_{y^E}^{y^E+h_y^E} dy,
  \end{gather*}
  where $\|\tilde{T} - \Pi_p \tilde{T} \|_{\Ltwo(\Domain_h)}^2 = \sum_{E \in \Domain_h} \|\tilde{T} - \Pi_p \tilde{T} \|_{\Ltwo(E)}^2$ was used.
  Consider a given desired element size at the boundary layer, denoted by $h_y^s$, defined relative to boundary layer width: $h_y^s = h_y \constlayerwidth$.
  Note that if one chooses $h_y^E$ so that
  \begin{gather}
    \left[\left(\frac{h_y^E}{\constlayerwidth}\right)^{p+1} \exp[-y^E / \constlayerwidth] \right]^2 \frac{1 - \exp[-2h_y^E / \constlayerwidth]}{2 h_y^E / \constlayerwidth} = [h_y^{p+1}]^2
    \label{eq:exponential-scaling}
  \end{gather}
  then one will have
  \begin{gather*}
    \|\tilde{T} - \Pi_p \tilde{T} \|_{\Ltwo(\Domain_h)}^2 \leq \sum_{E \in \Domain_h} [C h_y^{p+1}] ^2 \int_{y^E}^{y^E+h_y^E} dy = [C h_y^{p+1}]^2 L_y
    \quad \Rightarrow \quad \| \tilde{T} - \Pi_p \tilde{T} \|_{\Ltwo(\Domain_h)} \leq L_y C h_y^{p+1},
  \end{gather*}
  where $L_y$ denotes the length of the domain in the $y$ direction.
  Thus, scaling the elements according to \eqref{eq:exponential-scaling} results in the approximation bound retaining the desired decay behavior as $h_y$ goes to zero while allowing the mesh elements to grow exponentially in size away from the layer.
  Note that a similar scaling law will result by solving the optimization problem that minimizes the error bound over a set of element sizes.
  That all said, strictly enforcing \eqref{eq:exponential-scaling} requires solving the following transcendental equation for $h_y^E$:
  \begin{gather}
    (h_y^E)^{p+1} \exp[-y^E / \constlayerwidth] \left[\frac{1 - \exp[-2h_y^E / \constlayerwidth]}{2 h_y^E / \constlayerwidth}\right]^{1/2} = (h_y^s)^{p+1}.
    \label{eq:exponential-scaling-transcendental}
  \end{gather}
  Note that for $h_y^E / \constlayerwidth \ll 1$, enforcing \eqref{eq:exponential-scaling-transcendental} to leading order becomes
  \begin{gather}
    (h_y^E)^{p+1} \exp[-y^E / \constlayerwidth] = (h_y^s)^{p+1}
    \quad \Rightarrow \quad
    h_y^E = h_y^s \exp[y^E \sqrt{\kpar/\kperp}/(p+1)],
    \label{eq:exponential-scaling-inner}
  \end{gather}
  whereas enforcing \eqref{eq:exponential-scaling-transcendental} to leading order for $h_y^E / \constlayerwidth \gg 1$ becomes
  \begin{gather}
    \begin{aligned}
      (h_y^E)^{p+1} \exp[-y^E / \constlayerwidth] &[2 h_y^E / \constlayerwidth]^{-1/2} = (h_y^s)^{p+1}
      \quad \Rightarrow \quad \\
      &h_y^E = (h_y^s)^{\frac{p+1}{p+1/2}} \left[2\sqrt{\kpar/\kperp}\right]^{\frac{1}{2p+1}} \exp[y^E \sqrt{\kpar/\kperp}/(p+\tfrac{1}{2})].
    \end{aligned}
    \label{eq:exponential-scaling-outer}
  \end{gather}
  While the use of \eqref{eq:exponential-scaling-transcendental} directly or an approach that switches from \eqref{eq:exponential-scaling-inner} to \eqref{eq:exponential-scaling-outer} away from the layer would yield more efficiency, \eqref{eq:exponential-scaling-inner} is used throughout the mesh because it is easier to implement than \eqref{eq:exponential-scaling-transcendental} and is more stringent (i.e., has slower exponential growth) away from the layer than \eqref{eq:exponential-scaling-outer}.
  The exponentially refined mesh $\Domain_h$ that enforces \eqref{eq:exponential-scaling-inner} is constructed by first partitioning $\Domain$ into $m_x$ rectangles (i.e. $h_x=\pi/m_x$ and $h_y=1$), and then refining all rectangles $E$ in the $y$ direction into two smaller rectangles until
  \begin{gather}
    h_y^E \leq h_y^s \exp[y^E \sqrt{\kpar/\kperp} / (p+1)]
    \quad \text{and} \quad
    h_y^E \leq h_y^s \exp[(1-(y^E+h^E)) \sqrt{\kpar/\kperp} / (p+1)]
    \label{eq:exponential-uniform}
  \end{gather}

  \begin{figure}[tbp]
    \center
    \includegraphics[width=0.6\textwidth]{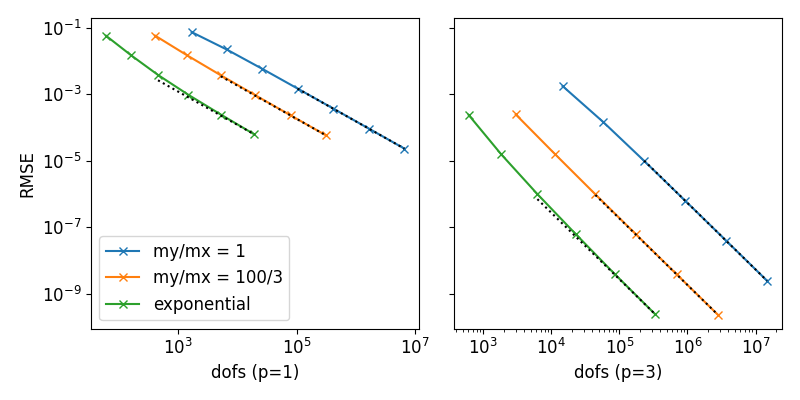}
    \caption{
      Efficiency diagram for the constant magnetic field problem \eqref{eq:toy} using uniform $m_x=m_y$ meshing ($h_x=h_y$), aspect ratio $m_y = 100/3 m_x$ ($h_x = 100h_y \, \pi/3$), and exponential meshing \eqref{eq:exponential-uniform} strategies (black dashed lines indicate theoretical power law).
      Note the aspect ratio and exponential refined meshes require orders of magnitude less degrees of freedom to achieve the same accuracy than is needed by the uniform meshing strategy.
    }
    \label{fig:variable-efficiency-study}
  \end{figure}

  The computational efficiency, measured as the error norm in \Cref{sec:verification} versus number of degrees of freedom (dofs), is obtained for uniform, aspect ratio, and exponential refined meshes with $\kpar/\kperp=10^4$.
  Both the aspect ratio and exponential refined meshes use $m_x = 3,6,12,24,48,96$ and $h_y^s = 1/ (10^2 \cdot m_x/3)$, with the aspect ratio meshes then using $h_y^E = h_y^s$ while the exponential refined meshes use \eqref{eq:exponential-uniform}.
  Note that one expects a theoretical power law dependence of error on the number of dofs $N$.
  Consider that $N$ is proportional to the number of elements in the domain, which is proportional to $1/(h_y^s)^2$ and thus $h_y^s \sim 1/\sqrt{N}$.
  Thus, one might expect the error in temperature to be proportional to $(h_y^s)^{p+1} \sim 1/N^{(p+1)/2}$ using the approximation theory bound \eqref{eq:interpolation-error} with $k=0$ and $l=1$.
  Figure \ref{fig:variable-efficiency-study} shows the efficiency results for piecewise linear ($p=1$) and piecewise cubic ($p=3$) function spaces.
  Both the aspect ratio and exponential refined meshes obtain the respective theoretical power law with far fewer dofs than required by the uniform meshes.
  Generally speaking, the aspect ratio meshes can attain the same accuracy solution as the uniform meshes with around an order of magnitude less dofs.
  Furthermore, the exponential refined meshes can obtain a solution of a given accuracy using an order of magnitude less dofs than the aspect ratio meshes and two to three orders of magnitude less dofs than the uniform meshes.

  \begin{figure}[tbp]
    \center
    \includegraphics[width=0.45\textwidth]{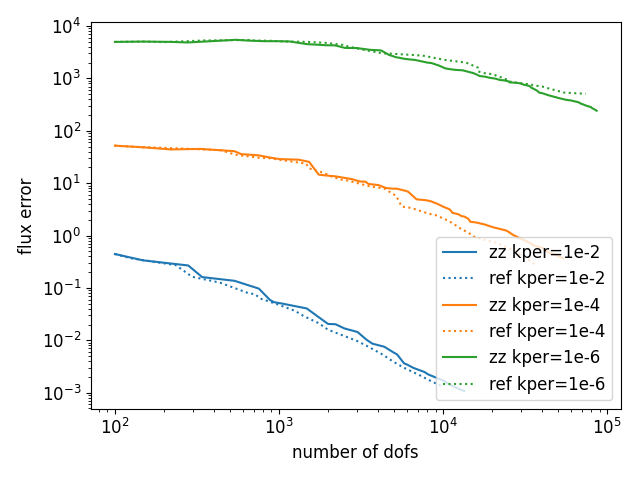}
    \caption{
      Efficiency diagram for threshold refinement using ``true'' flux error (reference solution) compared to using the ZZ estimator.
      The results show that the efficiency of AMR with threshold refinement using the ZZ estimator is similar to that using the reference flux error, and thus, the ZZ estimator is a sufficient error estimator for threshold refinement of our problem.
    }
    \label{fig:zz-vs-reference-efficiency}
  \end{figure}

\subsection{Mesh refinement for poloidally varying magnetic field problems}

  Note that while the aspect ratio and exponential refined meshes are very efficient for a constant magnetic field that is aligned with the mesh, these approaches do not immediately generalize to situations where either the mesh is not aligned with the magnetic field or when the magnetic field is dynamic in time.
  In particular, the aspect ratio mesh requires the mesh is aligned the magnetic field.
  The exponential refinement approach is more general in that it requires that the magnetic field is known \textit{a priori} but not that the mesh is aligned to the field.
  The caveat is that $y^E$ and $1 - (y^E+h^E)$ in \eqref{eq:exponential-uniform} must be replaced with some other measurement of the influence of the separatrix.
  The most general approach is adaptive mesh refinement (AMR), as it can be used to focus the refinement where it is needed regardless of whether the magnetic field is known \textit{a priori}.

\begin{figure}[tbp]
    \center
    \includegraphics[width=0.45\textwidth]{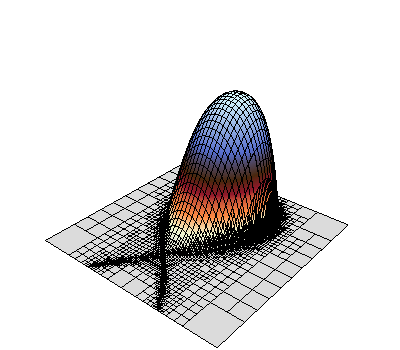}
    \includegraphics[width=0.45\textwidth]{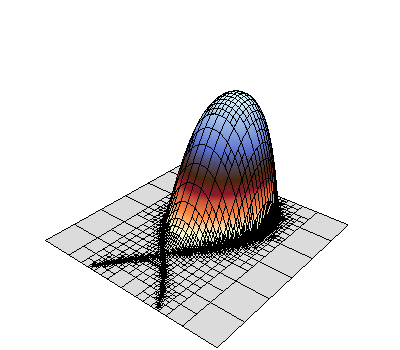}
    \caption{
      Solution profiles for $\kpar/\kperp = 10^4$ with AMR using the ZZ estimator (left) versus using reference flux error (right) for the single null problem.
      The  similarity between the AMR meshes indicate that the ZZ estimator is suitable for threshold refinement.
    }
    \label{fig:zz-vs-reference-solution}
  \end{figure}

  Recall from Section \ref{sec:methods-amr} that the ZZ refinement strategy combines a threshold strategy paired with the ZZ estimate of the flux error \eqref{eq:flux-error}.
  To first evaluate the performance of the ZZ error estimator independent of the performance of the threshold strategy, a reference solution is computed on a highly refined mesh for the single null problem \eqref{eq:two-wire} to provide the ``true'' flux error on each element to the element identification routine (in lieu of an analytic solution).
  The reference solutions for $\kpar/\kperp$ values of $10^2$, $10^4$, and $10^6$ are computed on uniform meshes consisting of $250^2$, $1000^2$ and $2000^2$ elements, respectively.
  Note that the reference mesh sizes are chosen so that (i) the reference mesh has elements that are smaller than the boundary layer and (ii) the results herein are qualitatively the same as when $125^2$, $500^2$, and $1000^2$ are used.
  The final piece of implementation detail is that the mesh refinement is halted when any element on the refined mesh becomes as small as the uniform elements in the reference solution.
  Figure \ref{fig:zz-vs-reference-efficiency} compares the efficiency of the AMR approach using the ZZ error estimate to the efficiency of the AMR approach using the reference flux error.
  Note that the efficiency of the AMR using the ZZ error estimate is qualitatively similar, and sometimes even better, than that of the AMR using the reference error.
  Figure \ref{fig:zz-vs-reference-solution} shows that the final meshes (and $\Temp$ field) are qualitatively the same for $\kpar/\kperp=10^4$ (i.e., the same regions are highly-refined).

  With the ZZ error estimator verified as a sufficient estimator for threshold AMR refinement, the efficiency of the approach is now measured against exponential and uniform refinement for the single null magnetic field \eqref{eq:two-wire}, double null magnetic field \eqref{eq:dbl-null}, and magnetic island \eqref{eq:mag-island} problems.
  Two different error metrics are considered for measuring the efficiency.
  The first is the error in temperature, defined as $\int_\Domain |\Temp - \Temp^*| \dVol$.
  The second is the error in flux, defined as before: $\int_\Domain \|\kappaTensor \nabla \Temp - \kappaTensor \nabla \Temp^*\|_2 \dVol$.
  Note that the errors are measured using the reference solution as $\Temp^*$ and quadrature rules to compute integrals.
  Recall that one might expect the error in temperature to be proportional to $h^{p+1} \sim 1/N^{(p+1)/2}$ using the approximation theory bound \eqref{eq:interpolation-error}.
  While a similar prediction for the flux error from \eqref{eq:interpolation-error} is not as direct, one still expects the error to be proportional to $h^p \sim 1/N^{p/2}$ following the same decrease in order of $h$ using $k=1$ in \eqref{eq:interpolation-error}.

\begin{figure}[tbp]
    \center
    \includegraphics[width=0.32\textwidth]{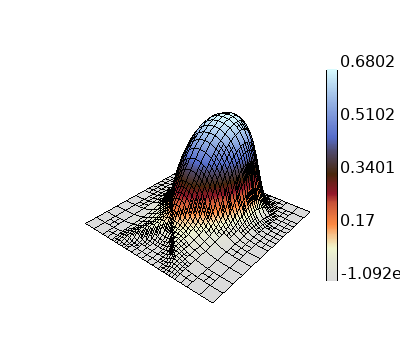}
    \includegraphics[width=0.32\textwidth]{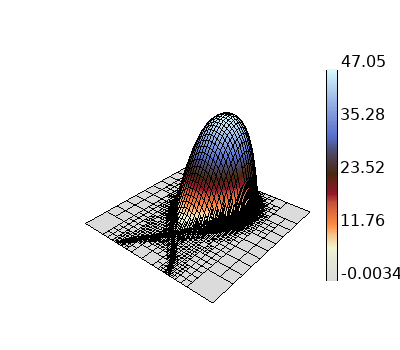}
    \includegraphics[width=0.32\textwidth]{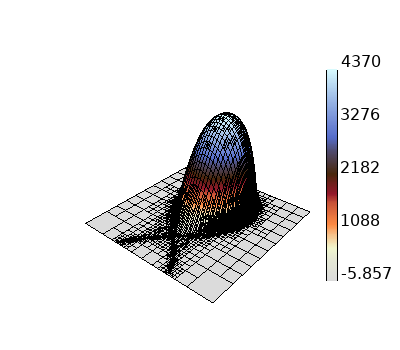}
    \\
    \includegraphics[width=0.32\textwidth]{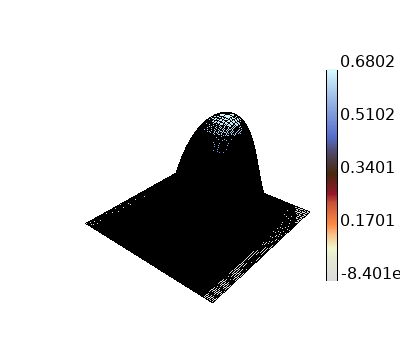}
    \includegraphics[width=0.32\textwidth]{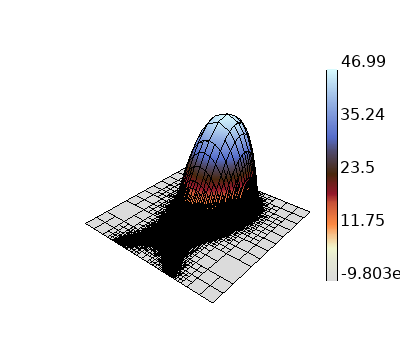}
    \includegraphics[width=0.32\textwidth]{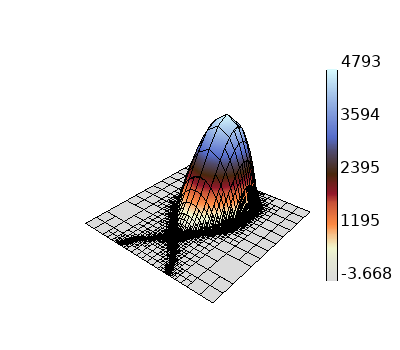}
    \caption{
      ZZ refinement (top) and exponential refinement (bottom) solutions for single null magnetic field problem with various $\kpar/\kperp$ values: $10^2$ (left), $10^4$ (center), $10^6$ (right).
    }
    \label{fig:two-wire-solutions}
  \end{figure}

  For the single null magnetic field problem \eqref{eq:two-wire}, the exponential refinement strategy is defined by first marking all elements that contain the separatrix that is implicitly defined by $A_z(\x)=A_z \big(\x_s)$, and then using the newly refined element size $h^s$ to mark all other elements where
  \begin{gather}
    h^E > h^s \exp[|A_z(\x_i) - A_z(\x_s)|\sqrt{\kpar/\kperp}/(p+1)]
    \label{eq:exponential-two-wire}
  \end{gather}
  with $\x_i$ being any of the quadrature points for the element $E$.
  Figure \ref{fig:two-wire-solutions} shows the ZZ and exponential refinement solutions with $\kpar/\kperp$ values of $10^2$, $10^4$, and $10^6$.
  Note that both solutions show temperature increasing proportional to the increase in $\kpar/\kperp$ and show refinement in the vicinity of the separatrix, although the vicinity is typically larger with exponential refinement.
  Another key difference is that the ZZ refinement approach refines more than the exponential refinement approach in the interior of the separatrix where the source drives the maximum temperature.
  Note this is consistent with the \eqref{eq:exponential-two-wire} being independent of the source in \eqref{eq:two-wire}.

\begin{figure}[tbp]
    \center
    \includegraphics[width=0.49\textwidth]{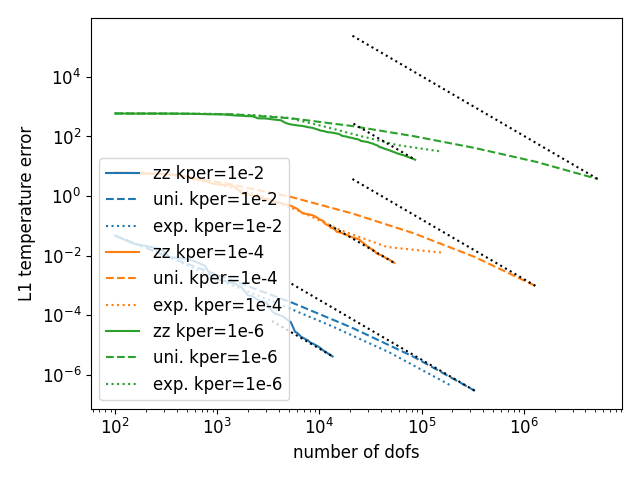}
    \includegraphics[width=0.49\textwidth]{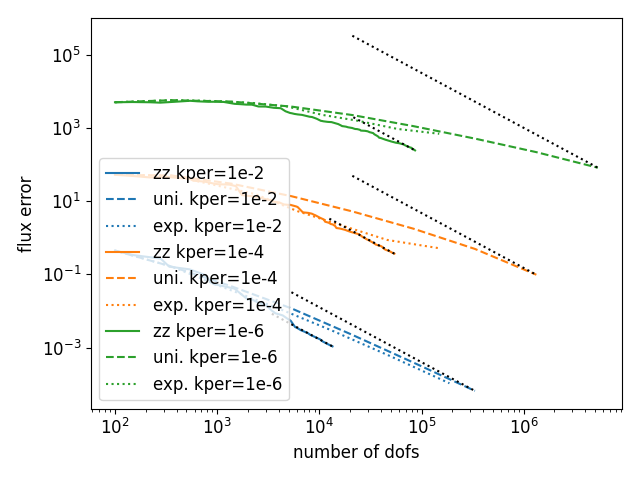}
    \caption{
      Efficiency diagrams for single null magnetic field problem with $\kpar/\kperp \in \{10^2, 10^4, 10^6\}$, using  $L_1$ temperature error (left) and sum of element flux errors defined in \eqref{eq:flux-error} (right).
      Results show that ZZ refinement tends to attain the theoretical power law (dotted black lines) with substantially fewer dofs than with uniform refinement, making it more efficient for large $\kpar/\kperp$.
    }
    \label{fig:two-wire-efficiency}
  \end{figure}

  Figure \ref{fig:two-wire-efficiency} investigates the efficiency of the ZZ and exponential refinement approaches, as well as that of uniform refinement for comparison.
  The efficiency results show a clear advantage for the ZZ refinement approach over that of uniform refinement, characterized by the ability to attain the theoretical power law dependence of error on dofs ($1/N^2$ and $1/N^{3/2}$, respectively) with significantly fewer dofs.
  Specifically, the power law dependence is attained with at least one or two orders of magnitude less dofs for $\kpar/\kperp$ values of $10^2$ and $10^4$, respectively.
  While neither solution attains the power law dependence for $\kpar/\kperp=10^6$ with the allotted number of dofs, the results indicate even more advantage at the higher $\kpar/\kperp$ ratio for the ZZ refinement solution.
  For the exponential refinement approach, a similar advantage in efficiency over uniform refinement is seen for anisotropy ratios large enough to limit the vicinity of the separatrix that is refined (i.e., $\kpar/\kperp$ values of $10^4$ and $10^6$).
  The efficiency advantage becomes degraded, however, once the dominant source of error becomes resolving the high temperatures inside the separatrix, which is something the continued refinement in the vicinity of the separatrix curve does not address.
  While the point at which the advantage is degraded likely depends on the choice of $|A_z(\x_i)-A_z(\x_s)|$ in \eqref{eq:exponential-two-wire}, there will still be a degradation point due to the independence of \eqref{eq:exponential-two-wire} on the source in the single null magnetic field problem \eqref{eq:two-wire}.

\begin{figure}[tbp]
    \center
    \includegraphics[width=0.32\textwidth]{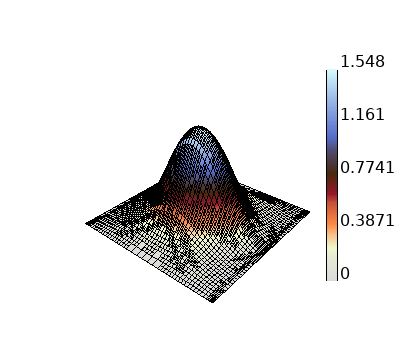}
    \includegraphics[width=0.32\textwidth]{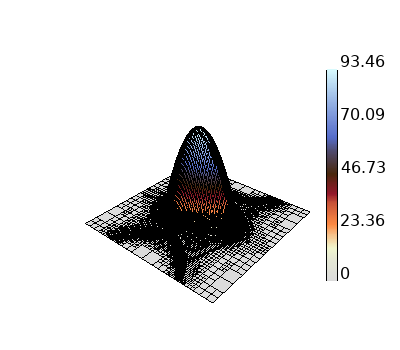}
    \includegraphics[width=0.32\textwidth]{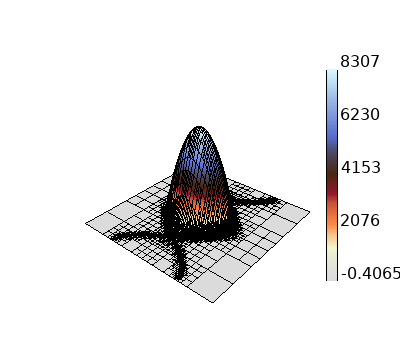}
    \\
    \includegraphics[width=0.32\textwidth]{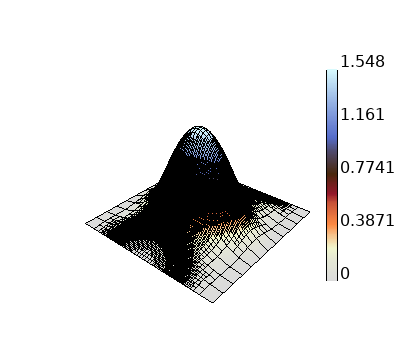}
    \includegraphics[width=0.32\textwidth]{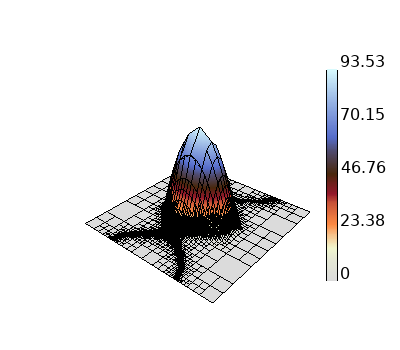}
    \includegraphics[width=0.32\textwidth]{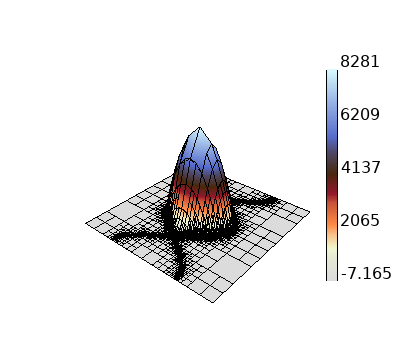}
    \caption{
      ZZ refinement (top row) and exponential refinement (bottom row) solutions for the double null magnetic field problem with $\kpar/\kperp$ values: $10^2$  (left),  $10^4$ (center),  $10^6$ (right).
    }
    \label{fig:dbl-null-solutions}
  \end{figure}

  For the double null magnetic field problem \eqref{eq:dbl-null}, the exponential refinement strategy is defined by again first marking all elements containing the separatrix defined by $A_z(\x) = A_z(\x_s)$, and then using the newly refined element size $h^s$ to mark all other elements where
  \begin{gather}
    h^E > h^s \exp \left[ \left| \frac{A_z(\x_i) - A_z(\x_s)}{A_z(\x_0) - A_z(\x_s)} \right| \sqrt{\kpar/\kperp}/(p+1)\right]
    \label{eq:exponential-mag-island}
  \end{gather}
  with $\x_i$ being any of the quadrature points for the element $E$.
  Figure \ref{fig:dbl-null-solutions} shows the ZZ and exponential refinement solutions for the double null magnetic field problem \eqref{eq:dbl-null} with $\kpar/\kperp$ values of $10^2$, $10^4$, and $10^6$.
  Like the single null magnetic field problem, both solutions show the temperature increasing proportional to $\kpar/\kperp$ and refinement in the vicinity of the separatrix.
  Also like the single null magnetic field problem, the ZZ refinement approach refines the interior of the separatrix where the source is driving the temperature whereas the exponential refinement approach is focused only on the vicinity of separatrix.
  While one might adjust the vicinity by choosing a different scaling and/or dependence of \eqref{eq:exponential-mag-island} on $A_z(\x)$, the independence of \eqref{eq:exponential-mag-island} on the source or boundary in the double null magnetic field problem \eqref{eq:dbl-null} means the approach can miss solution features that require additional resolution.

\begin{figure}[tbp]
    \center
    \includegraphics[width=0.49\textwidth]{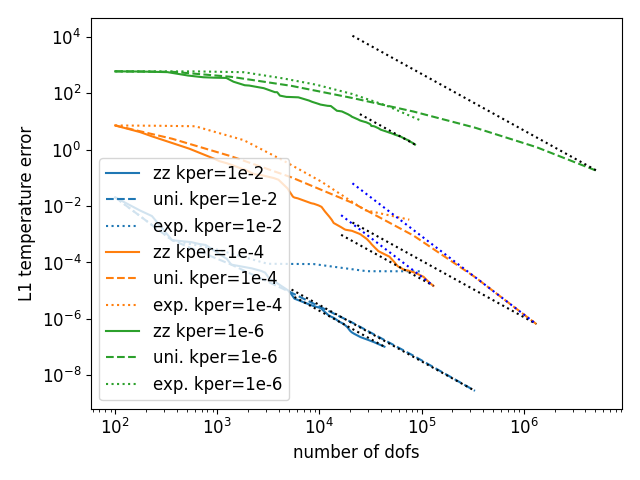}
    \includegraphics[width=0.49\textwidth]{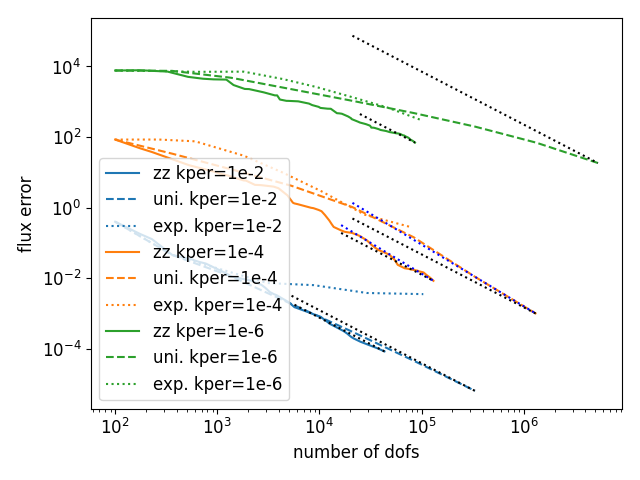}
    \caption{
      Efficiency diagrams for the double null magnetic field problem with various $\kpar/\kperp \in \{10^2, 10^4, 10^6\}$, using $L_1$ temperature error (left) and sum of element flux errors defined in \eqref{eq:flux-error} (right) compared to theoretical power law (dotted black lines) and empirical power law (dotted blue lines) when theoretical power law is exceeded.
      Results show that ZZ refinement starts to require substantially fewer dofs to attain the power law exhibited by uniform refinement once $\kpar/\kperp$ is large enough, indicating that ZZ refinement is more efficient for large $\kpar/\kperp$.
    }
    \label{fig:dbl-null-efficiency}
  \end{figure}

  Figure \ref{fig:dbl-null-efficiency} shows the corresponding efficiency results for the double null magnetic field problem.
  For $\kpar/\kperp = 10^2$, there is no significant advantage for the ZZ or exponential refinement approaches over the uniform refinement approaches as all approaches attain the theoretical power law dependence around the same small number of dofs.
  Note that the ZZ refinement solution for $\kpar/\kperp=10^2$ in Figure \ref{fig:dbl-null-solutions} looks close to a uniform mesh, which is consistent with the ZZ refinement efficiency more or less matching that of uniform refinement.
  Note also that the exponential refinement approach again sees degraded efficiency once the dominant error comes from a region in the domain that is not included in the vicinity of the separatrix.
  For $\kpar/\kperp = 10^4$, the ZZ and uniform refinement approaches both attain a power law dependence that is stronger than the theoretical predictions for both error metrics, with the power of $1/N$ exceeding $2$ and $3/2$, respectively.
  That said, the ZZ refinement approach does attain the same power law dependence shown by the uniform refinement approach using only around half the dofs.
  The exponential refinement approach does not exhibit a sustained power law dependence and is the least efficient of the three approaches, likely due to not sufficiently refining the regions away from the separatrix.
  As with the single null problem, no approach clearly attains the power law dependence for $\kpar/\kperp=10^6$ with the allotted number of dofs, although the results indicate a strong advantage of the ZZ refinement approach at the higher $\kpar/\kperp$ ratio.

\begin{figure}[tbp]
    \center
    \includegraphics[width=0.32\textwidth]{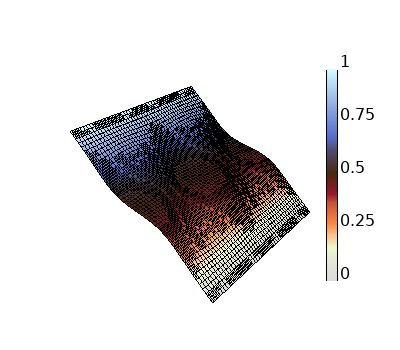}
    \includegraphics[width=0.32\textwidth]{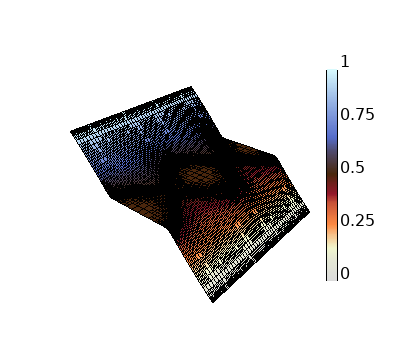}
    \includegraphics[width=0.32\textwidth]{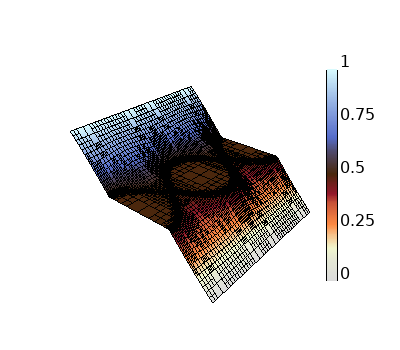}
    \\
    \includegraphics[width=0.32\textwidth]{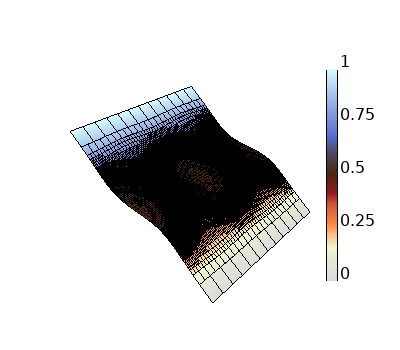}
    \includegraphics[width=0.32\textwidth]{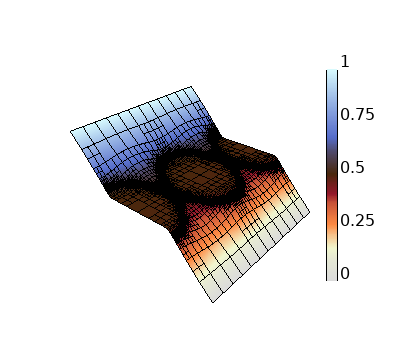}
    \includegraphics[width=0.32\textwidth]{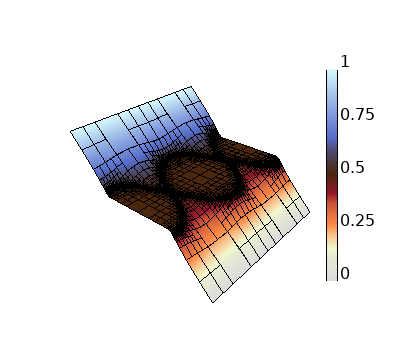}
    \caption{
      ZZ refinement (top row) and exponential refinement (bottom row) solutions for the magnetic island problem with various $\kpar/\kperp$ values: $10^2$ (left), $10^4$ (center), $10^6$ (right).
    }
    \label{fig:mag-island-solutions}
  \end{figure}

  For the magnetic island problem \eqref{eq:mag-island}, the exponential refinement strategy is the same as for the double null magnetic field problem, i.e., using \eqref{eq:exponential-mag-island}.
  Figure \ref{fig:mag-island-solutions} shows the ZZ and exponential refinement solutions for $\kpar/\kperp$ values of $10^2$, $10^4$, and $10^6$.
  Unlike the single null and double null magnetic field problems, the magnetic island problem does not have a source inside the separatrix.
  As such, higher values of $\kpar/\kperp$ do not result in higher temperatures but, instead, in a more constant temperature problem inside the separatrix.
  Again, the key difference between the ZZ and exponential refinement approaches is that the former focuses refinement both along the separatrix and in regions of large temperature gradients whereas the latter focuses on the separatrix due to \eqref{eq:exponential-mag-island} being unable to account for source or boundary effects.

  Figure \ref{fig:mag-island-solutions} shows the corresponding efficiency results for the magnetic island problem.
  Like the double null magnetic field problem, the ZZ and uniform refinement approaches are about as efficient as the other for $\kpar/\kperp$, due to ZZ refinement being essentially uniform (see Figure \ref{fig:mag-island-solutions}), with the ZZ refinement approach gaining a substantial advantage for larger $\kpar/\kperp$.
  As with all the poloidally varying magnetic field problems, the advantage can again be explained as the ZZ refinement approach requiring substantially less dofs to obtain the same power law dependence as the uniform refinement approach eventually does.
  The exponential refinement approach again shows a degradation in efficiency, likely due to the lack of refinement outside the separatrix where there is a large temperature gradient.
  As with all the poloidally varying magnetic field problems, the exponential refinement lacks a consistent efficiency advantage due to its inability to refine solution features that result from sources and/or boundary conditions.
  \begin{figure}[tbp]
    \center
    \includegraphics[width=0.49\textwidth]{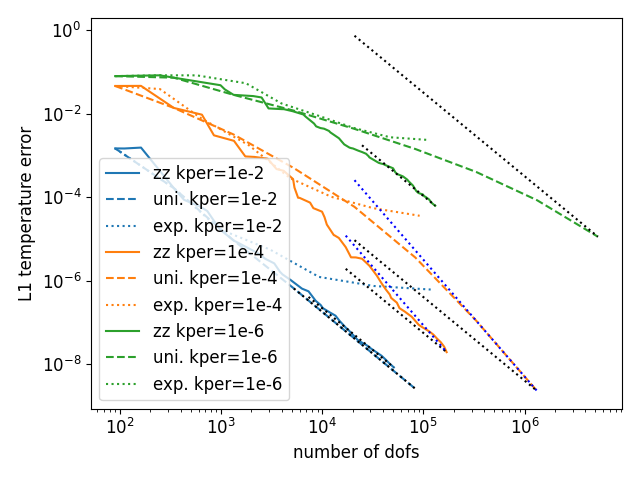}
    \includegraphics[width=0.49\textwidth]{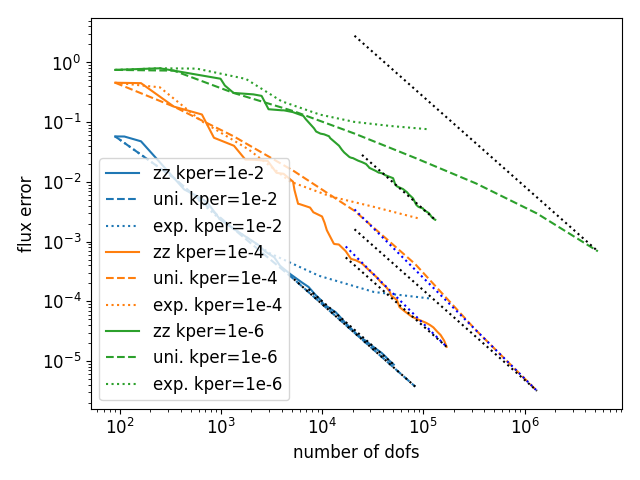}
    \caption{
      Efficiency diagrams for the magnetic island problem with $\kpar/\kperp \in \{10^2, 10^4, 10^6\}$, using $L_1$ temperature error (left) and flux error (right) compared to theoretical power law (dotted black lines) and empirical power law (dotted blue lines) once it exceeds the theoretical power law.
      Results show that, once $\kpar/\kperp$ is large enough, ZZ refinement requires substantially fewer dofs to attain the power law exhibited by uniform refinement, indicating that ZZ refinement is more efficient for large $\kpar/\kperp$.
    }
    \label{fig:mag-island-efficiency}
  \end{figure}

%% file: sections/discussion.tex
\def\Nelements{\elements}
\def\Jmax{J}
\def\jindex{j}

 \def\xrad{y}

 \def\hrad{ h_r}
 \def\hpol{h_p}
 \def\htor{h_t}
 \def\hsurf{h_s}
 \def\hpar{h_\|}

 \def\nrad{n_r}
 \def\npol{n_p}
 \def\ntor{n_t}
 \def\nsurf{n_s}
 \def\npar{n_\|}

\def\Cost{{\mathcal C}} %Computational Cost

 \Cref{sec:results} demonstrated how adaptive mesh refinement can substantially improve the computational efficiency in solving anisotropic diffusion problems, such as the linear steady-state problem \eqref{eq:ssp}.
 Those results, however, were limited to a handful of idealized poloidally varying magnetic field problems.
 As such, some discussion is warranted on how much improvement in computation efficiency one can expect for more general magnetic fields, particularly with large anisotropy ratios.
 First, theoretical scaling laws for the degrees of freedom requirements with mesh refinement and general magnetic fields are presented.
 Then, the dependence of the matrix condition number on the anisotropy ratio is discussed with a focus on the implications on the cost of solving the underlying discrete anisotropic diffusion problem.
 The resulting scaling laws are then combined to get an approximation of the total cost $\Cost$ under uniform, isotropic, and anisotropic adaptive mesh refinement.

\subsection{Scaling of dofs for constant magnetic field}

  To develop mesh refinement scaling laws for general magnetic fields, first consider the constant magnetic field $\bhat = \xhat$ that causes a boundary layer at $\xrad=0$ with width $\wrad$.
  MFEM supports a variety of mesh refinement approaches, including to uniformly refine a quadrilateral/hexahedral element (i) isotropically into four/eight quadrilateral/hexahedral elements or (ii) anisotropically in two/four quadrilateral/hexahedral elements.

  Consider applying either of these uniform refinement approach to a single element mesh, where the elements containing the boundary layer are refined along the $\xrad$ direction.
  If the refinement process is repeated $\Jmax$ times, the resulting mesh will have nodes at $\xrad=0$ and at $\xrad_\jindex = L/2^\jindex$ for $\jindex=0,\ldots,\Jmax$.
  Choose $\Jmax$ to be the smallest value such that the boundary layer is resolved by at least one element, or equivalently that $\xrad_\Jmax \leq \wrad$ and $y_{\Jmax-1} > \wrad$.
  Note that $\Jmax$ is the smallest integer that is greater than $\log(\Len/\wrad)/\log(2)$.
  The quantity of interest is how many elements $\Nelements$ result from $\Jmax \approx \log(\Len/\wrad)/\log(2)$ iterations of refinement, which can be written as the sum, $\Nelements = \Nelements_1 + \ldots + \Nelements_\Jmax$, where $\Nelements_\jindex$ is the number of elements from the $\jindex$-th round of refinement (i.e., elements that have nodes at $y_\jindex$ and $y_{\jindex-1}$).

  First, consider anisotropic refinement in two dimensions (2D), i.e., where a quadrilateral element is uniformly refined into two smaller quadrilaterals.
  After the first refinement, $\jindex=1$, the mesh now has two elements, one of which contains the boundary and one which does not contain the boundary layer.
  Note that the element that does not contain the boundary layer will not be further refined: $\Nelements_1 = 1$.
  The element that contains the boundary layer is refined into two elements: one that contains the boundary layer and one that does not.
  Note that the new element that does not contain the boundary layer will not be further refined: $\Nelements_2 = 1$.
  Generalizing, one can see that $\Nelements_\jindex= 1$ for $\jindex < \Jmax$.
  At the last round of refinement, $\jindex=\Jmax$,  all the new elements are no longer refined: $\Nelements_\Jmax = 2$.
  Thus, the number of elements for the two-dimensional isotropic refinement mesh is
  \begin{gather}
    \Nelements_{\anisoamr} = \Jmax+1 \approx \frac{ \log (\Len/\wrad) }{\log(2)} + 1 \sim \log \sqrtkparkperp
    \label{eq:anisotropic-element-scaling}
  \end{gather}
  for $\kpar/\kperp \gg 1$.
  Note that this result also applies to three dimensions (3D).

  Now, consider isotropic refinement in two dimensions (2D),  i.e., where a quadrilateral element is uniformly refined into four smaller quadrilaterals.
  After  the first refinement level, $\jindex=1$, the mesh now has four elements, two of which contain the boundary and two which do not contain the boundary layer.
  Note that the two elements that do not contain the boundary layer will not be further refined: $\Nelements_1 = 2$.
  The two elements that contain the boundary layer are each refined into four elements: two that contain the boundary layer and two that do not.
  Note that the new elements that do not contain the boundary layer will not be further refined: $\Nelements_2 = 2 \cdot \Nelements_1 = 4$.
  The four elements that contain the boundary layer are each refined into four elements: two that contain the boundary layer and two that do not.
  Again note that the new elements that do not contain the boundary layer will not be further refined: $\Nelements_3 = 2 \cdot \Nelements_2$.
  Generalizing, one can see that $\Nelements_\jindex = 2 \cdot \Nelements_{\jindex-1} $ for $\jindex < \Jmax$.
  At the final round of refinement, $\jindex=\Jmax$, all the elements are no longer refined: $\Nelements_\Jmax =  4 \cdot \Nelements_{\Jmax-1} $.
  Thus, the number of elements for the two-dimensional isotropic refinement mesh is
  \begin{gather}
    \Nelements_{\isoamr}^{2D} = \sum_{j=1}^{J} 2^j + 2^{J} = 3 \cdot 2^\Jmax - 2 \approx 3 (\Len/\wrad) - 2 \sim \sqrtkparkperp
    \label{eq:isotropic-element-scaling-2d}
  \end{gather}
  for $\kpar / \kperp \gg 1$.
  Note that in three dimensions, where each hexahedron is refined into eight smaller hexahedrons, one has $\Nelements_1 = 4$, $\Nelements_\jindex = 4 \cdot \Nelements_{\jindex-1} $ for $\jindex < \Jmax$, and $\Nelements_\Jmax =  8 \cdot \Nelements_{\Jmax-1}$.
  Thus, the number of elements for the three-dimensional isotropic refinement mesh is
  \begin{gather}
    \Nelements_{\isoamr}^{3D} = \sum_{\jindex=1}^\Jmax 4^\jindex + 4^\Jmax = \tfrac{7}{3} (2^\Jmax)^2 - \tfrac{4}{3} \approx \tfrac{7}{3} (L/w)^2 - \tfrac{4}{3} \sim \kpar/\kperp.
    \label{eq:isotropic-element-scaling-3d}
  \end{gather}
  Thus, one can clearly see the advantage of using anisotropic refinement with field aligned meshes for large $\kpar / \kperp$, especially in three-dimensional simulations.

\subsection{Scaling of dofs for general geometry}

  The expression for the interpolation error \eqref{eq:interpolation-error} implies that in order to restrict the maximum error, it is beneficial to exponentially pack the mesh near the boundary layer.
  For a simple Cartesian mesh, the mesh points would lie at the point $\xrad_j = \Len e^{-j/\nrad}$, where $\Len$ is the overall mesh size and the  parameter $\nrad$ controls the rate at which cells are refined.
  For field aligned meshes, the previous section found, $\nrad=1/ \log(2)$, but this rate could potentially differ for more general geometric situations.
  The number of  elements in the radial direction required to reach the boundary layer width $\wrad$ is
  \begin{align}
   \Jmax = \nrad \log{(\Len/\wrad)}\sim \nrad   \log \sqrtkparkperp.
  \end{align}

   %\emph{Isotropic Refinement.}
 The lesson learned from our analysis is that, for isotropic refinement one also needs to pack at the smallest radial distance $\wrad$ in   the directions tangent to flux surfaces.
 In this case, refining a factor of 2 at each level yields $2^{d-1}$ more elements per level for a total of $\sim 2^{J(d-1)}=(L/w)^{(d-1)}$ elements.
 This leads to the estimate
  \begin{align}
    \begin{aligned}
      \Nelements_{\isoamr} &\sim (\Len/\wrad)^{\dim-1 }
      \sim (\sqrtkparkperp)^{(\dim-1)}
      \\
      \Nelements_{\isoamr}/\Nelements_{\uni}&\sim  \wrad/\Len
      \sim (\kappa_\perp/\kappa_\|)^{1/2}
    \end{aligned}\label{eq:Nelements_iso-amr}
  \end{align}
which generalizes the results of the previous section.

\iffalse
  For the 2D problems under consideration, the number of elements is
  \begin{align} \label{eq:Nelements_iso-amr_2D}
    \Nelements_{\isoamr}^{2D} &\sim (\Len/  \wrad)   \sim \sqrtkparkperp
  \end{align}
 and, in 3D, the number of elements required is even larger
  \begin{align} \label{eq:Nelements_iso-amr_3D}
    \Nelements_{\isoamr}^{3D} &=(\Len /\wrad)^2 \sim    (\kappa_\|/\kappa_\perp)  .
  \end{align}
\fi

  %\emph{Anisotropic Refinement.}
  In principle, anisotropic refinement can perform much better than isotropic, if the elements can become field aligned at the smallest scales, because one does not need to resolve such small length scales tangent to the boundary layer.  Assume that this still requires $\nsurf=\Len/\hsurf$ points per dimension tangent to boundary layer, where $\hsurf$ is the distance needed to resolve the surface.
  This yields the estimate
  \begin{align}
    \begin{aligned}
      \Nelements_{\anisoamr} &\sim  (\Len/\hsurf)^{\dim-1} \nrad \log(\Len/\wrad) \sim \nrad  \nsurf^{\dim-1}  \log\sqrtkparkperp
      \\
      \Nelements_{\anisoamr} /\Nelements_{\uni} &\sim   (\wrad^d/\hsurf^{\dim-1}\Len)\nrad \log(\Len/\wrad)
      \sim  \nrad  \nsurf^{2\dim-1} (\kappa_\perp/\kappa_\|)^{\dim/2} \log\sqrtkparkperp.
    \end{aligned}
  \end{align}
  In 2D, $\npol=\Len/\hpol$ where $\hpol$ is the smallest distance in poloidal direction.  This is set by the needs of accuracy and is related to the poloidal curvature of the surface, $\aminor$, rather than the boundary layer width.  Thus, the estimate is simply
  \begin{align} \label{eq:Nelements_aniso-amr_2D}
    \Nelements_{\anisoamr} &\sim (\Len / \hpol )  \nrad \log(\Len/\wrad) \sim \nrad\npol  \log\sqrtkparkperp
 %   \\
%    \Nelements_{\anisoamr} /\Nelements_{\uni}&\sim (\wrad^2 / \hpol \Len)  \nrad \log(\Len/\wrad)     \sim  \nrad\npol (\kappa_\perp/\kappa_\|)  \log\sqrtkparkperp
.
  \end{align}
  In 3D, the number of points required in the toroidal direction is $ \ntor=\Len/\htor$ where $\htor$ is the smallest distance required in the toroidal direction.
  Again, this is set by the need to resolve the toroidal curvature, $\Rmajor$, rather than by the boundary layer width.
  Thus, the result is
  \begin{align} \label{eq:Nelements_aniso-amr_3D}
    \Nelements_{\anisoamr} &\sim  ( \Len^2/ \hpol\htor ) \nrad \log(\Len/\wrad)
     \sim \nrad \npol\ntor \log \sqrtkparkperp
%    \\
%    \Nelements_{\anisoamr} /\Nelements_{\uni}&    \sim     ( \wrad^3/ \hpol\htor\Len ) \nrad \log(\Len/\wrad)  \sim \nrad\npol\ntor  (  \kappa_\perp/\kappa_\|)^{3/2}   \log\sqrtkparkperp
.
  \end{align}

\begin{figure}[tbp]
  \centering
    \includegraphics[width=0.49\textwidth]{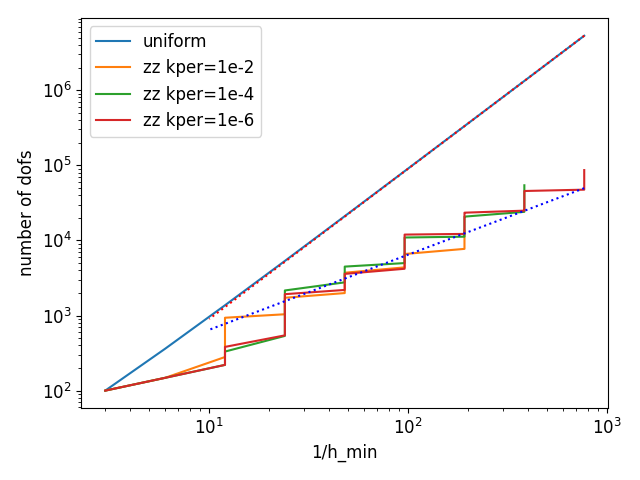}
    \includegraphics[width=0.49\textwidth]{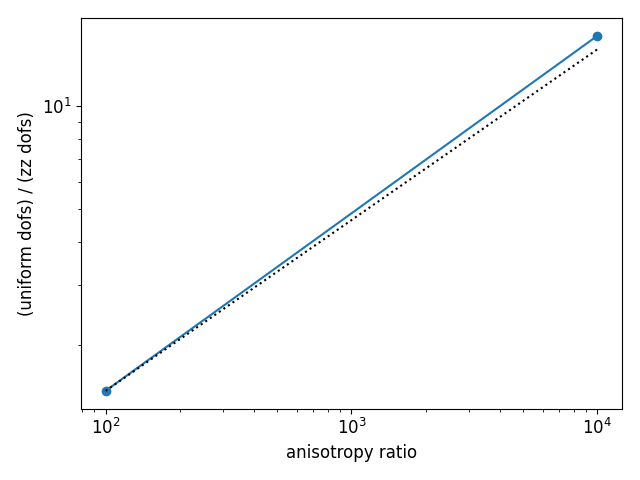}
    \caption{
      Scaling of dofs versus the reciprocal of minimum element size $1/h_{\rm min}$ for each iteration of refinement (left) and versus the magnetic anisotropy ratio $\kpar/\kperp$ (right) for the single null problem.
      The dofs versus $1/h_{\rm min}$ scaling follows the power law $(\kpar/\kperp)^2$ (dotted red line) for uniform refinement and follows the power law $\kpar/\kperp$ (dotted blue line) for ZZ refinement.
      The dofs ratio of uniform refinement to ZZ refinement versus $\kpar/\kperp$ scaling follows the predicted power law $\sqrtkparkperp$ (dotted black line).
    }
    \label{fig:two-wire-dofs}
  \end{figure}

  %Comparison to ZZ estimator performance
  To compare against the theoretical scaling laws, Figure \ref{fig:two-wire-dofs} shows the number of dofs, which are proportional to the number of elements, for the single null problem \eqref{eq:two-wire}.
  The minimum element size $h_{\rm min}$ is noted as the meshes for $\kpar/\kperp \in \{10^2, 10^4, 10^6\}$ are refined according to the uniform and ZZ refinement strategies.
  Note that when the number of dofs at each iteration of refinement are compared to $1/h_{\rm min}$, the uniform and ZZ refinement follow power laws of $(\kpar/\kperp)^2$ and $(\kpar/\kperp)$, respectively.
  The number of dofs for both uniform and ZZ refinement once $h_{\rm min}$ becomes smaller than $\sqrt{\kperp/\kpar}$ is also noted for each value of $\kpar/\kperp \in \{10^2, 10^4\}$.
  Note that the ratio of those dofs scales as $\sqrtkparkperp$, as predicted by \eqref{eq:Nelements_aniso-amr_2D}.
  Specifically, the ZZ refinement mesh resolves the layer with a factor of almost $2$ fewer dofs for $\kpar/\kperp = 10^2$ and a factor of almost $20$ fewer dofs for $\kpar/\kperp=10^4$.
  Extrapolating these results to an anisotropy ratio of $10^6$, one would predict a reduction in dofs by a factor of around $200$.

\subsection{Scaling of dofs for filamentary structures}
  Although this work is focused on steady-state transport, in the future it would be useful to extend the results to the consideration of dynamic plasma turbulence and MHD activity.
  It is well known that it is advantageous for the instabilities responsible for driving plasma turbulence and MHD preferentially form elongated filamentary structures along field lines.
  This is perhaps best exemplified by edge localized modes that are active in the pedestal of a tokamak operating in the high-performance (H-mode) confinement regime.

 Isotropic AMR is an excellent strategy for capturing low-dimensional structures in higher dimensional volumes.
  In order to resolve a small region with perpendicular width $\wrad$ at the corner of a $d$-dimensional hypercube, one requires $\Jmax$ iterations with $\Nelements \sim 2^\dim$ elements for each iteration for a total number of elements $\sim  2^\dim \Jmax $.
  If the point is in the interior of the domain, this requires $\sim 2^{2\dim} \Jmax\sim2^{2\dim}\log{(\Len/\wrad}) $ elements.
   %Exact result to resolve a corner for 2x refinement at each level is [J (2^d-1) + 1] for a total of 2^d[J (2^d-1) + 1]

  For a filamentary structure, there is only one direction along the filament.
  Thus, isotropic AMR in the perpendicular directions is a good strategy for resolving the core of the filament.
  For an anisotropic AMR procedure that can align with the filament at the smallest scales, then one can assume that the number of elements required along the field line is $\npar=\Len/\hpar$, where $\hpar$ is determined by the parallel field line curvature.
  This yields
  \begin{align}
  \Nelements_\anisoamr  &\sim (\Len/\hpar)  2^{2(d-1)}  \log(\Len/\wrad) \sim  \npar 2^{2(d-1)}  \log \sqrtkparkperp
  %\\
 % \Nelements_\anisoamr/\Nelements_\uni  &\sim (\wrad^{d}  / \hpar \Len^{\dim-1} ) \log(\Len/\wrad) \sim \npar  (\kperp/\kpar)^{\dim/2}   \log \sqrtkparkperp .
  \end{align}
For isotropic AMR, the direction along the filament must also be resolved to width $\wrad$.
Each iteration layer has $\sim 2^{2(\dim-1)}$ perpendicular regions with $\sim 3\times 2^\Jmax$ elements required to cover the parallel direction, for a total of $3\cdot 2^{2(d-1)}(\Len/\wrad)$ elements.
%Exact result is [ (2^d-1) (2^J - 2) + 2^{J+1} 2^d ] ~ 3x2^{d+J}
Hence, the result is
  \begin{align}
  \Nelements_\isoamr  &\sim   3\cdot 2^{2(d-1)} (\Len/\wrad) \sim  3\cdot 2^{2(d-1)} \sqrtkparkperp
  %\\
 % \Nelements_\isoamr/\Nelements_\uni  &\sim  (4\wrad  /\Len)^{\dim-1} \sim    (4\kperp/\kpar)^{(\dim-1)/2}   .
  \end{align}
Thus, isotropic AMR performs relatively well in this scenario.

\subsection{Solver Iterations}

  With an iterative approach to solving the underlying sparse linear system that represents the discrete steady state problem \eqref{eq:ssp-fem}, the total cost of the solution (i.e., number of flops) is on the order of the product of the number of iterations required and the total number of degrees of freedom (dofs).
  Thus, while the number of dofs is a significant factor in the cost of solving the problem, the dependence of the number of required iterations on the condition number of the matrix in \eqref{eq:ssp-fem} is another important factor.
  For example, the convergence rate of the conjugate gradient method depends on the reciprocal of the condition number.
  To discuss the matrix condition number in \eqref{eq:ssp-fem}, consider that the eigenfunction problem for the constant magnetic field problem \eqref{eq:toy} has solutions of $\exp[i 2 q_x x]\exp[i 2\pi q_y y]$, for non-zero integers $q_x$ and $q_y$, with corresponding eigenvalues $\lambda_{q_x q_y} = \kpar (2q_x)^2 + \kperp (2\pi q_y)^2$.
  While the eigenfunctions and eigenvalues of the finite-dimensional operator in \eqref{eq:ssp-fem} will not be the same as those for the operator in \eqref{eq:toy}, the difference between the two should be proportional to the mesh spacing.
  Thus, let $q_{x,h}$ and $q_{y,h}$ represent the largest frequencies supported by $\Domain_h$ so that the matrix condition number in \eqref{eq:ssp-fem} is approximated by
  \begin{gather}
    \frac{\max_{q_x,q_y} |\lambda_{q_x q_y}|}{\min_{q_x,q_y} |\lambda_{q_x q_y}|} = \frac{\kpar (2q_{x,h})^2 + \kperp (2\pi q_{y,h})^2}{\kpar 2^2 + \kperp (2\pi)^2} =  \frac{(\kpar/\kperp) (2q_{x,h})^2 + (2\pi q_{y,h})^2}{(\kpar/\kperp)2^2 + (2\pi)^2}.
    \label{eq:condition}
  \end{gather}

\begin{figure}[tbp]
    \center
    \includegraphics[width=\textwidth]{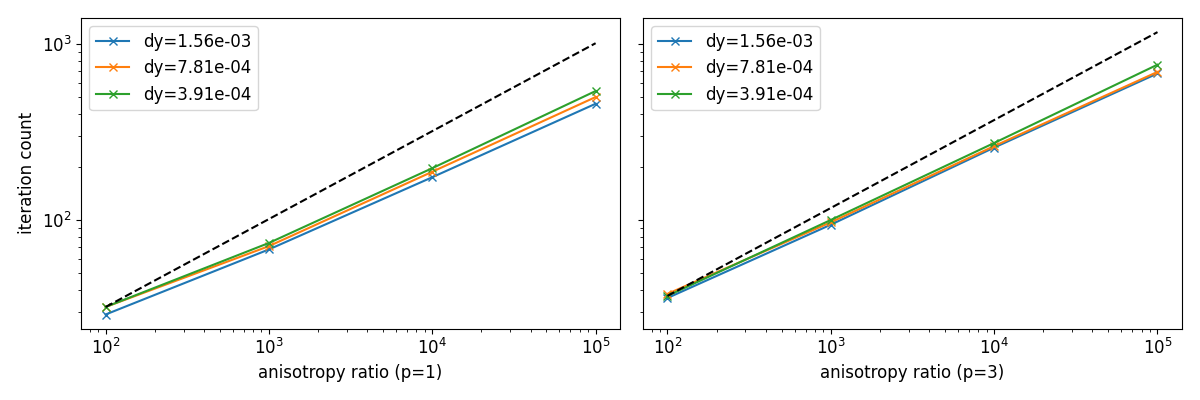}
    \caption{
      Dependence of iteration count on anisotropy ratio $\kpar/\kperp$ for uniform refinement meshes as the boundary layer in the constant magnetic field problem \eqref{eq:toy} is resolved.
      Results show a power law dependence that is bounded above by the theoretical estimate of $\sqrtkparkperp$ shown as a dashed black line.
    }
    \label{fig:iterations-constant-field-uniform}
  \end{figure}

  Consider when $\Domain_h$ consists of $\Mgrid \times \Mgrid$ uniform elements, resulting in $q_{x,h} = q_{y,h}$.
  Denote that largest supported frequency by $q_h$ and note that \eqref{eq:condition} is now
  \begin{gather}
    \frac{\max_{q_x,q_y} |\lambda_{q_x q_y}|}{\min_{q_x,q_y} |\lambda_{q_x q_y}|} = q_h^2
    \label{eq:condition-uniform}
  \end{gather}
  Recall that Section \ref{sec:results} found that the boundary layer of width $\sqrt{\kperp/\kpar}$ must be resolved by the elements in $\Domain_h$ to reach the asymptotic regime, which would result in $q_h = \sqrtkparkperp$.
  Therefore, it is important to note that as larger values of $\kpar/\kperp$ are considered, the number of conjugate gradient iterations required to obtain the solution on a uniform mesh is expected to grow proportional to $\sqrtkparkperp$ (i.e., the square-root of the reciprocal of the matrix condition number).
  To test this $\sqrtkparkperp$ scaling, the iteration counts for conjugate gradient preconditioned with BoomerAMG are collected for \eqref{eq:toy}.
  Figure \ref{fig:iterations-constant-field-uniform} shows that the iteration counts do indeed exhibit a power law dependence on anisotropy ratio bounded above by the $\sqrtkparkperp$ scaling law for both linear ($\poly=1$) and cubic ($\poly=3$) function spaces.

\begin{figure}
    \center
    \includegraphics[width=\textwidth]{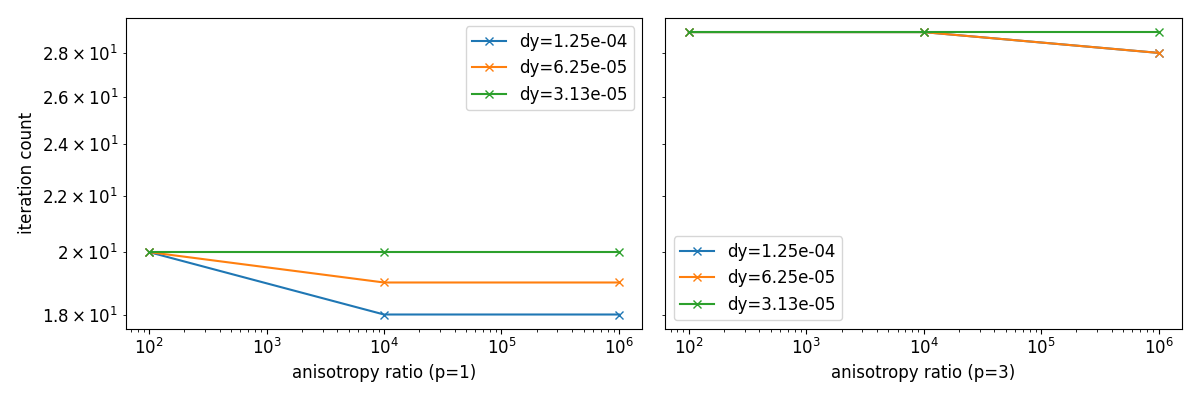}  \\
    \caption{
      Dependence of iteration count on anisotropy ratio $\kpar/\kperp$ for ratio refinement meshes as the boundary layer in the constant magnetic field problem \eqref{eq:toy} is resolved.
      Results show near independence of iteration count on anisotropy ratio.
    }
    \label{fig:iterations-constant-field-ratio}
  \end{figure}

  When the mesh is aligned with the magnetic field, the ratio and exponential refinement approaches investigated in Section \ref{sec:results} can potentially avoid iteration count growth with anisotropy ratio.
  Consider a mesh where $m_x = \alpha m_y$ and $m_y=\sqrtkparkperp$.
  The largest supported frequencies on $\Domain_h$ are now $q_{x,h} = \alpha \sqrtkparkperp$ and $q_{y,h} = \sqrtkparkperp$.
  Thus, the matrix condition number is now
  \begin{gather}
    \frac{\max_{q_x,q_y} |\lambda_{q_x q_y}|}{\min_{q_x,q_y} |\lambda_{q_x q_y}|} = \frac{\alpha^2 (\kpar/\kperp)^2 2^2 + (\kpar/\kperp) (2\pi)^2}{(\kpar/\kperp) 2^2 + (2\pi)^2},
    \label{eq:condition-ratio}
  \end{gather}
  which is $\mathcal{O}(1)$ if $\alpha \sim \sqrt{\kperp/\kpar}$ for $\kpar/\kperp \gg 1$.
  One would then expect the iteration count on such a ratio refinement mesh, including the $m_y/m_x = 100/3$ mesh for $\kpar/\kperp = 10^4$, to be almost independent of growth in anisotropy ratio.
  To test the independence, the iteration counts for conjugate gradient preconditioned with BoomerAMG are again collected for \eqref{eq:toy}.
  Figure \ref{fig:iterations-constant-field-ratio} shows that indeed the iteration counts for both linear ($\poly=1$) and cubic ($\poly=3$) function spaces do indeed appear to be independent of $\kpar/\kperp$.

  The condition number estimates in \eqref{eq:condition-uniform} and \eqref{eq:condition-ratio} can be generalized beyond field aligned meshes for constant magnetic fields.
  Generally speaking, the condition number of the discrete parallel and perpendicular Laplacian operators scale as $\kpar (\Len_\| / h_\|)^2$ and $(\dim-1) \kperp (\Len_\perp / h_\perp)^2$, respectively, where $\dim$ is the dimension, $(\Len_\|, h_\|)$ are the longest and shortest wavelengths present in the parallel direction, and $(\Len_\perp,h_\perp)$ are the longest and shortest wavelengths present in the perpendicular direction.
  Assume the mesh resolves the boundary layer with width $\wrad$ (i.e., $h_\perp = \wrad$) and that the perpendicular directions interact with the domain boundary so that $\Len_\perp$ is approximately the domain length (i.e., $\Len_\perp \approx \Len$).
  The condition number of the combined anisotropic Laplacian for such a problem is approximately
  \begin{gather}
    \frac{\kpar (1/h_\|)^2 + (\dim-1) \kperp (1/\wrad)^2}{\kpar (1/\Len_\|)^2 + (\dim-1) \kperp (1/\Len)^2}.
    \label{eq:condition-general}
  \end{gather}
  Consider first problems where the parallel directions also interact with the boundary so that $\Len_\| = \Len$.
  The condition number estimate \eqref{eq:condition-general} is now
  \begin{gather*}
     (\Len/\wrad)^2 \frac{\kpar (\wrad/h_\|)^2 + (\dim-1) \kperp}{\kpar + (\dim-1) \kperp}
     =
     \frac{\kpar}{\kperp} \frac{\kpar (\wrad/h_\|)^2 + (\dim-1) \kperp}{\kpar + (\dim-1) \kperp}.
  \end{gather*}
  Note if a uniform mesh is used where $h_\| \approx h_\perp = \wrad$, the condition number scales approximately as $\kpar/\kperp$, which is consistent with \eqref{eq:condition-uniform}.
  Note if a ratio refined mesh is used where $h_\| \approx \sqrtkparkperp \wrad$, the condition number scales approximately independent of $\kpar/\kperp$, which is consistent with \eqref{eq:condition-ratio}.

  Consider now problems where the parallel directions do not interact with the boundary so that $\Len_\| = \infty$.
  The condition number estimate \eqref{eq:condition-general} is now
  \begin{gather*}
    (L/\wrad)^2 \frac{\kpar (\wrad/h_\|)^2 + (\dim-1) \kperp }{(\dim-1) \kperp}
    =
    \frac{\kpar}{\kperp} \frac{\kpar (\wrad/h_\|)^2 + (\dim-1) \kperp }{(\dim-1) \kperp}.
  \end{gather*}
  Note that if a uniform mesh is now used where $h_\| \approx \wrad$, the condition number scales approximately as $(\kpar/\kperp)^2$.
  Thankfully, it is observed that the single null \eqref{eq:two-wire}, double-null \eqref{eq:dbl-null}, and magnetic island \eqref{eq:mag-island} problems all exhibit iteration counts that scale as $\sqrtkparkperp$.
  This behavior is seen for both uniform refinement and ZZ refinement meshes, which is welcome news as even though ratio refinement theoretically can restore the condition number scaling to $\kpar/\kperp$ when the parallel directions do not interact with the boundary, a ratio refined mesh approach for poloidally varying magnetic fields may be quite nontrivial.
  Even once one has the field-aligned mesh, the dofs themselves may not be perfectly field-aligned even in an element with edges aligned.

\begin{figure}[tbp]
    \center
    \includegraphics[width=\textwidth]{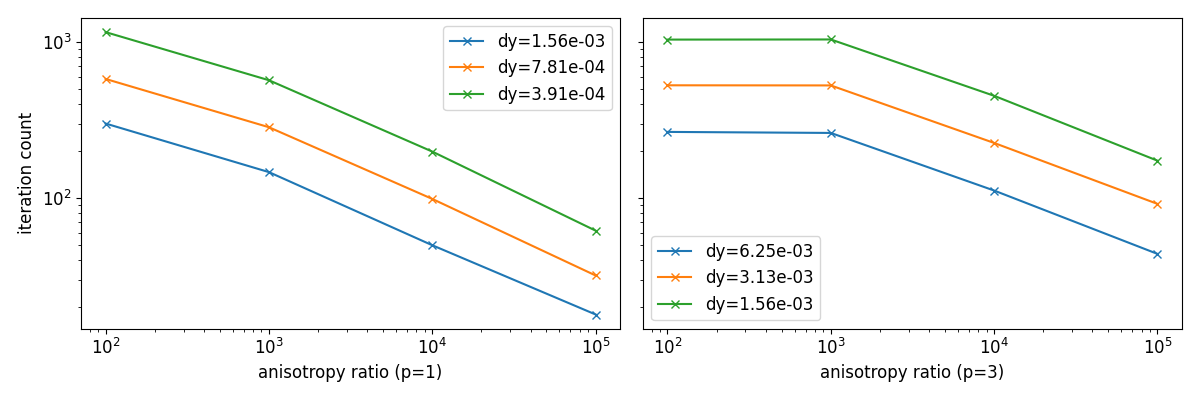}
    \caption{
      Dependence of iteration count on anisotropy ratio $\kpar/\kperp$ for uniform refinement meshes when ILU preconditioning is used for the  constant magnetic field problem \eqref{eq:toy}.
      Results show a substantial decrease in iteration count over BoomerAMG preconditioning, with the counts actually decreasing with larger anisotropy ratio.
    }
    \label{fig:iterations-constant-field-uniform-ilu}
  \end{figure}

  A successful solver approach for poloidally varying magnetic fields will therefore likely require both meshing and preconditioning strategies.
  As an example, ILU preconditioning results in iteration counts that do not grow with increasing magnetic anisotropy even on uniform meshes.
  In certain cases, including when applied to the constant magnetic field problem \eqref{eq:toy}, the iteration count actually decreases with increasing magnetic anisotropy (see Figure \ref{fig:iterations-constant-field-uniform-ilu}).
  The downsides are that ILU preconditioning is not generally scalable across multiple MPI ranks and that the same benefit did not appear to extend to non field-aligned meshes in our initial testing.

\subsection{Computational Cost}

  To summarize, one can estimate the total cost to solution by multiplying the scaling laws for dofs by the empirically measured scaling law for iteration count.
  For uniformly refined meshes, the cost would be
  \begin{equation}
    \Cost_\uni \sim (\Len/\wrad)^{d+1} \sim (\sqrtkparkperp)^{d+1}.
  \end{equation}
  For isotropic adaptive refined meshes, the cost is improved by the factor $\Len/\wrad \sim \sqrtkparkperp$
  \begin{equation}
    \Cost_\isoamr \sim (\Len/\wrad)^{d} \sim (\sqrtkparkperp)^d .
  \end{equation}
  For anisotropic adaptive refined meshes, the cost could potentially be as low as
  \begin{equation}
    \Cost_\anisoamr \sim (\Len/\wrad) \log (\Len/\wrad) \sim \sqrtkparkperp \log \sqrtkparkperp .
  \end{equation}

%% file: sections/conclusion.tex
This work investigated the efficiency of various mesh refinement strategies for boundary layers formed by the extreme magnetic anisotropy found in tokamak reactors.
While verifying the implementation of a finite element discretization in MFEM against an analytic solution for a test problem, we noted the theoretical convergence rate is not attained for any polynomial order until the mesh elements containing the boundary layer are small enough to capture the layer (see Figure \ref{fig:convergence-test}).
We found such behavior consistent with established approximation theory that indicates the error in the polynomial finite element solution will depend on the polynomial order and mesh size in such a way that the mesh size must be small enough in the layer to counter the increasing magnitude of higher-order derivatives of the exact solution.
The same approximation theory result indicates that substantially less mesh refinement is needed in the direction along the layer relative to across the layer.
Furthermore, the approximation theory result indicates less mesh refinement is needed for elements far away from the layer than for elements within the layer.
As such, we utilized two variable mesh strategies for meshes that are aligned with the magnetic field: one that used uniform rectangular elements with an aspect ratio of the layer width, and a second that chooses the aspect ratio of the rectangular elements to capture the exponential nature of the approximation theory bound.
We found that the uniform aspect ratio mesh attains the same accuracy as the uniform rectangular mesh using an order of magnitude less degrees of freedom (dofs) and that the exponential refined rectangular mesh attains another order of magnitude reduction in the number of dofs needed (see Figure \ref{fig:variable-efficiency-study}).

The variable mesh refinement approaches do not directly generalize to magnetic fields that do not align with the mesh, such as the spatially (and temporally) varying fields in a tokamak reactor.
As an example, the exponential refinement approach in Figure \ref{fig:variable-efficiency-study} can be used for magnetic fields with a given potential function after the $y_E$ input is replaced with a function of the potential value at $y_E$.
  \add{Note that the use of such a function based on the magnetic flux function is the current standard refinement approach in the field.}
For a more general refinement approach, we explored an adaptive mesh refinement (AMR) approach than refines elements that exhibit an estimated error larger than a specified threshold, with the estimated error obtained by the Zienkiewicz and Zhu approach.
After empirically verifying that the Zienkiewicz and Zhu error estimation leads to similar refinement behavior as error estimation using a highly refined reference solution, we measured the efficiency of the AMR and exponential refinement approaches against uniform refinement for three magnetic fields that model the spatial variation of poloidal fields expected in a tokamak.
For the single-null tokamak magnetic field geometry, we found the AMR approach attained the same accuracy as the uniform mesh using at least an order of magnitude less dofs (see Figure \ref{fig:two-wire-efficiency}), with the efficiency gain increasing with stronger anisotropy.
The exponential refinement approach also sees increased computational efficiency for higher anisotropy ratios while solution error along the separatrix remains dominant.
Once the solution error resolving the interior of the separatrix becomes dominant, the solution accuracy stagnates resulting in the loss of the efficiency advantage.
For the double null magnetic field that models the poloidal field in a double-null tokamak reactor, we found the AMR approach had at least an order of magnitude gain in efficiency once the anisotropy was strong enough whereas the exponential refinement approach did not exhibit a substantial improvement in efficiency (see Figure \ref{fig:dbl-null-efficiency}).
The magnetic island field that models a perturbation to the poloidal field showed similar results as with the double null magnetic field: significant efficiency gains for the AMR approach and no significant efficiency gain for the exponential refinement approach (see Figure \ref{fig:mag-island-efficiency}).
For all three fields, the efficiency advantage for the AMR approach is characterized by the AMR solution requiring less dofs than
  \change{the uniform solution}
  {either the uniform or exponential solutions}
before attaining the power law dependence of error on dofs exhibited by the uniform solution.
  \add{Given that the AMR approach is both more efficient for these fields at higher anisotropy ratios than the exponential strategy and does not require locating the separatrices and defining a distance function based on the magnetic flux function, we recommend the AMR approach be considered instead of the exponential strategy when simulating more complex magnetic geometries.}

While we explored anisotropy ratios up to $10^6$, the anisotropy in a tokamak reactor is expected to be much stronger.
We found that extending our experiments to anisotropy ratios beyond $10^6$ quickly becomes infeasible because of the increasing number of dofs required to resolve the layer compounded by the number of iterations required by the conjugate gradient solver used with algebraic multigrid preconditioning.
We derived scaling laws that showed the number of dofs for isotropic adaptive mesh refinement is less than uniform refinement by a factor of the square-root of the anisotropy ratio \eqref{eq:Nelements_iso-amr}, whereas the number of dofs for anisotropic adaptive mesh refinement is even smaller with only a logarithmic dependence on the square-root of the anisotropy ratio for both two-dimensional \eqref{eq:Nelements_aniso-amr_2D} and three-dimensional \eqref{eq:Nelements_aniso-amr_3D} problems.
For the field-aligned test problem used to verify the finite element discretization, we derived a theoretical estimate that the condition number associated with uniform refinement grows as the anisotropy ratio \eqref{eq:condition-uniform} and empirically verified the resulting prediction that the conjugate gradient iteration count, with algebraic multigrid preconditioner, grows as the square-root of the anisotropy ratio result (see Figure \ref{fig:iterations-constant-field-uniform}).
We similarly derived a theoretical estimate that the condition number associated with ratio refinement, where the rectangular elements have an aspect ratio equal to the layer width, is independent of the anisotropy ratio \eqref{eq:condition-ratio} and empirically verified that the iteration count is independent of anisotropy ratio (see Figure \ref{fig:iterations-constant-field-ratio}).
Such independence of iteration count and anisotropy ratio is crucial to generalize beyond ratio refined meshes, such as those produced by AMR, that can resolve spatially and temporally varying magnetic fields.
We note that ILU preconditioning does indeed result in independence of and substantial reduction in iteration count on a uniform mesh for the field-aligned test problem (see Figure \ref{fig:iterations-constant-field-uniform-ilu}); however, ILU preconditioning does not readily scale across multiple MPI ranks to tackle the larger number of dofs required by more realistic problems.
Thus, in the future, it would be desirable to explore other promising solver techniques, e.g. that use line smoothing or geometric multigrid along with the direct solution of the resulting coarsened problems.